\documentclass[reqno]{article}
\usepackage{eucal,amsmath,amsfonts,amsthm}
\usepackage[all]{xy}

\usepackage{vmargin}
\setpapersize{USletter}
\setmarginsrb{2.5cm}{1.75cm}{2.5cm}{1.75cm}%
{12pt}{12pt}{12pt}{24pt}

\newcommand{\spn}{\operatorname{span}}
\newcommand{\Adm}{\operatorname{Adm}}
\newcommand{\ov}{\overline}
\def\lo{\rightarrow}
\newcommand{\Pc}{{\mathcal{P}}}
\newcommand{\calo}{{\mathcal{O}}}
\def\cnabla{\check{\nabla}}
\newcommand{\pai}[2]{\langle #1, #2\rangle}
\newcommand{\setdef}[2]{\{#1 \; | \; #2\}}

\newcommand{\cinfty}[1]{C^\infty(#1)}
\newcommand{\map}[3]{#1\colon#2\rightarrow#3}
\newcommand{\set}[2]{\left\{\,#1\left.\vphantom{#1#2}\,\right\vert\,#2\,
                \right\}}

\newcommand{\real}{\mathbb{R}}
\newcommand{\integer}{\mathbb{Z}}
\renewcommand{\natural}{\mathbb{N}}

\newcommand{\invclos}[1]{\operatorname{\overline{Lie}}(#1)}
\newcommand{\symclos}[1]{\operatorname{\overline{Sym}}(#1)}
\newcommand{\symcloss}[2]{\operatorname{\overline{Sym}}^{#1}(#2)}
\newcommand{\grad}{\operatorname{grad}}
\newcommand{\psinv}{\operatorname{pseudoinv}}

\newcommand{\br}[2]{[#1,#2]} 
\newcommand{\symprod}[2]{\langle #1: #2\rangle}
\newcommand{\pd}[2]{\frac{\partial #1}{\partial #2}}
\newcommand{\at}[1]{\Big|_{#1}}
\newcommand{\vectorfields}[1]{\mathfrak{X}(#1)}
\newcommand{\Sec}[1]{\operatorname{Sec}(#1)}
\let\sec\Sec
\newcommand{\Hc}{{\mathcal{H}}}
\newcommand{\Gc}{{\mathcal{G}}}
\newcommand{\Hor}{\operatorname{Hor}}
\newcommand{\Ver}{\operatorname{Ver}}
\newcommand{\hor}{\operatorname{hor}}
\newcommand{\ver}{\operatorname{ver}}
\newcommand{\Br}{\operatorname{Br}}
\newcommand{\famY}{\mathcal{Y}}
\newcommand{\famB}{\mathcal{B}}
\newcommand{\famX}{\mathcal{X}}
\newcommand{\Id}[1]{\operatorname{Id}_{#1}}
\newcommand{\prol}[2][\,]{\CMcal{T}_{#1}#2}
\newcommand{\prolmap}[1]{\CMcal{T} #1}
\newcommand{\V}{{\CMcal{V}}}
\newcommand{\X}{{\CMcal{X}}}
\newcommand{\sode}{{\textsc{sode}}}
\newcommand{\spV}{^{\scriptscriptstyle V}}

\def\prb{\prol[]{\tau}}
\newcommand{\CC}{\Gamma} 
\newcommand{\SC}{f} 
\newcommand{\D}{{\mathcal{D}}}

\newcommand{\C}{{\mathcal{C}}}
\newcommand{\R}{{\mathcal{R}}}
\newcommand{\ada}[1]{\operatorname{ad}^\dag_{#1}}
\newcommand{\G}{\mathfrak{g}}

\renewcommand{\SC}{S} 

\newtheorem{contar}{}[section]
\newtheorem{definition}[contar]{Definition}
\newtheorem{proposition}[contar]{Proposition}
\newtheorem{theorem}[contar]{Theorem}

\newtheorem{lemma}[contar]{Lemma}
\newtheorem{remark}[contar]{Remark}

\begin{document}

\title{Mechanical control systems on Lie algebroids}

\author{Jorge Cort\'es\\
  Coordinated Science Laboratory\\
  University of Illinois at Urbana-Champaign\\
  1308 W Main St, Urbana, IL 61801, USA\\
  \texttt{jcortes@uiuc.edu}
  \and Eduardo Mart{\'\i}nez\\
  Departamento de Matem\'atica Aplicada\\
  Universidad de Zaragoza\\
  50009 Zaragoza, Spain\\
  \texttt{emf@unizar.es}}

\date{\today}
\maketitle

\begin{abstract}
  This paper considers control systems defined on Lie algebroids.
  After deriving basic controllability tests for general control
  systems, we specialize our discussion to the class of mechanical
  control systems on Lie algebroids. This class of systems includes
  mechanical systems subject to holonomic and nonholonomic
  constraints, mechanical systems with symmetry and mechanical systems
  evolving on semidirect products. We introduce the notions of linear
  connection, symmetric product and geodesically invariant subbundle
  on a Lie algebroid.  We present appropriate tests for various
  notions of accessibility and controllability, and analyze the
  relation between the controllability properties of control systems
  related by a morphism of Lie algebroids.
\end{abstract}

\section{Introduction}

One of the basic problems in control theory is that of deciding the
local controllability properties of a given system. Roughly speaking,
the local controllability problem consists of finding appropriate
conditions guaranteeing that the set of reachable states starting from
an initial point is open, i.e., that the system can move locally in
\emph{any} direction. Deciding the controllability properties of a
system is an \emph{a priori} question that one needs to have addressed
before being able to undertake other control problems such as motion
planning and trajectory generation.  The controllability problem has
received a great deal of attention during the last decades
(see~\cite{HH:82,HN-AJvdS:90,EDS:98,HJS:87} and references therein).
In particular, researchers have undertaken a thorough study of control
systems with a rich geometric structure such as mechanical and
homogeneous systems, and made use of their special properties to
accomplish accurate modeling settings and sharp analysis results.
Specific class of control problems include simple mechanical
systems~\cite{ADL-RMM:95c}, systems subject to nonholonomic
constraints~\cite{AMB-MR-NHMC:92,FB-MZ:01d,ADL:97a}, systems invariant
under the action of a Lie group of
symmetries~\cite{JC-SM-JPO-HZ:02,SDK-RMM:95,SMM-PEC:84,SM-JC:03},
systems enjoying special homogeneity
properties~\cite{JC-SM-FB:01m,MK:95,PAV-JWB:03}, systems evolving on
semidirect products~\cite{JS:02}, and more.

One of the features which imposes a separate study for each class of
systems is the lack of a unified framework.  For instance, it is well
known that a Lagrangian system invariant under the action of a Lie
group of symmetries can be reduced to the quotient space induced by
the action, but the reduced dynamics is not the one that corresponds
to the reduced Lagrangian. Recent investigations have lead to a
unifying geometric framework to overcome this drawback.  It is
precisely the underlying structure of Lie algebroid on the phase space
what allows a unified treatment. This idea was introduced by
Weinstein~\cite{AW:96} in order to define a Lagrangian formalism which
is general enough to account for the various types of systems. A
symplectic formalism was later introduced for Lagrangian~\cite{EM:01a}
and Hamiltonian systems~\cite{EM:01b}.  One of the advantages of the
Lie algebroid formalism is the possibility of establishing appropriate
maps (morphisms) between two systems that respect the structure of the
phase space, and allow to relate their respective control properties.
The underlying idea (in a similar way to the notion of
\emph{abstraction}~\cite{GJP-GL-SSS:00}) is that the property of
interest will be easier to decide for one of the systems, and that by
means of the morphism one will be able to infer the same knowledge for
the other system.

In this paper we study the controllability problem for systems affine
in the inputs evolving on a bundle $E$ with an underlying Lie
algebroid structure, $\tau: E \rightarrow M$.  We build on previous
work on local controllability of general systems~\cite{HJS:87} and of
mechanical systems~\cite{ADL-RMM:95c} to derive tests to check the
accessibility and controllability properties of control systems
evolving on a Lie algebroid. Throughout the paper, we pay special
attention to what we term \emph{mechanical control systems on a Lie
  algebroid}. This class of systems embraces a variety of different
situations that can occur when analyzing mechanical systems, such as
the ones mentioned above.  Building on the notion of prolongation
$\prol{E} \rightarrow E$ of the Lie algebroid $E \rightarrow M$
introduced in~\cite{EM:01a}, we develop all the necessary differential
geometric tools enabling an intrinsic treatment of the second-order
dynamics associated with mechanical systems.  We focus on the set of
reachable points in the base manifold $M$ and in the bundle space $E$
starting from states which belong to the zero section of $E
\rightarrow M$. We make use of the geometric homogeneity properties of
the controlled equations, which turn out to greatly simplify the
accessibility and controllability computations.  We carefully describe
the relation between the controllability properties of control systems
that are related by a morphism of Lie algebroids.  As a result of the
generality of the approach, we are able to present in a unified way
previous work on the configuration accessibility and controllability
properties of simple mechanical control
systems~\cite{FB-NEL-ADL:00,JC-SM-JPO-HZ:02,ADL-RMM:95c,SM:02} (see
also~\cite{FB-ADL:04} for a comprehensive overview).  Regarding
systems evolving on semidirect products, the application of the Lie
algebroid approach renders novel tests which are valid in slightly
more general settings than the ones considered in~\cite{JS:02}.  We
also extend notions such as fiber controllability to what we call
\emph{controllability with regards to a manifold} and develop
conditions to check this property.

In the course of the preparation of this manuscript, we came across
the recent research effort~\cite{PAV-JWB:03}. This reference, which is
close in spirit to this work, analyzes the controllability properties
of so-called ``1-homogeneous control systems'' evolving on a vector
bundle.  However, it deals with vector fields with values in the
tangent bundle of the vector bundle, as opposed to deal with the
formalism of Lie algebroids and their prolongations.  This choice of a
higher-dimensional phase space makes necessary to resort to additional
geometric tools such as Ehresmann connections in order to describe the
structure of the accessibility algebra. We think that the Lie
algebroid approach accommodates the same level of generality, while
enabling in general a more concise treatment of the controlled
dynamics.

The paper is organized as follows. In Section~\ref{se:lie-algebroids}
we present some basic facts on Lie algebroids. We also discuss in
detail the notion of linear connection on a Lie algebroid, including
the generalization of the Levi-Civita connection and the constrained
connection.  In Section~\ref{se:prolongation}, we introduce the
prolongation of a Lie algebroid and develop the differential geometry
of horizontal sections, homogeneity, \sode\ sections and geodesically
invariant subbundles. In Section~\ref{se:general-control-systems} we
study nonlinear affine control systems whose drift and input vector
fields are associated with some sections of a Lie algebroid.  This
apparent restriction is not such, since most physical systems can be
casted into this form. We formulate the conditions for local
accessibility and controllability in terms of Lie brackets of
sections, and we study the effect of a morphism of Lie algebroids on
these properties. In Section~\ref{se:mechanical-control-systems} we
introduce the class of mechanical control systems defined on a Lie
algebroid. We show that the notion of affine connection control system
can be generalized to the setting of Lie algebroids, thus providing a
general framework to study the controllability properties of these
systems. We introduce the notions of local base controllability and
controllability with regard to a manifold, and we obtain computable
sufficient conditions to check them. We also study the effect of
morphisms of Lie algebroids in simplifying the controllability
analysis.  These results are later applied in
Section~\ref{se:applications} to simple mechanical systems defined on
a manifold, simple mechanical systems with symmetry and systems
defined on semidirect products and orbits of group actions.
Section~\ref{se:conclusions} presents some concluding remarks. We have
gathered in an appendix some basic notions for control systems defined
on manifolds. A final remark is that the summation convention over
repeated indexes is understood throughout the paper.

\section{Preliminaries on Lie algebroids}\label{se:lie-algebroids}
\setcounter{equation}{0}

In this section we introduce some known notions and develop new
concepts concerning Lie algebroids that are necessary for the further
developments. We refer the reader to~\cite{AC-AW:99,KM:87} for
thorough studies of Lie groupoids, Lie algebroids and their role in
differential geometry.  Let $M$ be an $n$-dimensional manifold and let
$\map{\tau}{E}{M}$ be a vector bundle with $\ell$-dimensional fibers.
A structure of \emph{Lie algebroid on $E$} is given by a Lie algebra
structure on the $\cinfty{M}$-module of sections of the bundle,
$(\Sec{E},\br{\cdot\,}{\cdot})$, together with a homomorphism
$\map{\rho}{E}{TM}$ of vector bundles (called the \emph{anchor map})
satisfying the compatibility condition
\[
\br{\sigma_1}{F\sigma_2} = F\br{\sigma_1}{\sigma_2} + \bigl(
\rho(\sigma_1)F \bigr) \sigma_2 \, .
\]
Here $F$ is a smooth function on $M$, $\sigma_1$, $\sigma_2$ are
sections of $E$ and we have denoted by $\rho(\sigma)$ the vector field
on $M$ given by $\rho(\sigma)(m)=\rho(\sigma(m))$.  The homomorphism
$\rho$ is called the \emph{anchor map}.  From the compatibility
condition and the Jacobi identity, it follows that the map
$\sigma\mapsto\rho(\sigma)$ is a Lie algebra homomorphism from
$\Sec{E}$ to $\vectorfields{M}$.

It is convenient to think of a Lie algebroid $\map{\tau}{E}{M}$ as a
substitute of the tangent bundle of $M$. In this way, one regards an
element $a$ of $E$ as a generalized velocity, and the actual velocity
$v$ is obtained when applying the anchor to $a$, i.e., $v=\rho(a)$.

The image of the anchor map, $\rho(E)$, defines an integrable smooth
generalized distribution on $M$.  Therefore, $M$ is foliated by the
integral leaves of $\rho(E)$, which are called the \emph{leaves of the
  Lie algebroid}.  A curve $a: [t_0,t_1] \rightarrow E$ is said to be
\emph{admissible} if $\dot{m}(t)=\rho(a(t))$, where $m(t)=\tau(a(t))$,
$t \in [t_0,t_1]$.  It follows that $a(t)$ is admissible if and only
if the curve $m(t)$ lies on a leaf of the Lie algebroid, and that two
points are in the same leaf if and only if they are connected by (the
base curve of) an admissible curve.

A Lie algebroid is said to be \emph{transitive} if it has only one
leaf, which is obviously equal to $M$. It is easy to see that $E$ is
transitive if and only if $\rho$ is surjective. If $E$ is not
transitive, then the restriction of the Lie algebroid to a leaf
$L\subset M$, $E_{|L}\lo L$ is transitive. In the latter case, one can
show that $E_{|L}$, and hence $\ker{\rho}$, has constant dimension.
We will say that a Lie algebroid is \emph{locally transitive at a
  point $m \in M$} if $\rho_m: E_m \rightarrow T_mM$ is surjective. In
this way, $m$ is contained in a leaf of maximal dimension.

Given a local basis $\{e_\alpha\}_{\alpha=1}^\ell$ of sections of $E$
defined on an open set $V \subset M$, we can write $a = y^\alpha
e_\alpha (\tau(a))$ for any $a \in E$ such that $\tau (a) \in V$. If
$(x^i)$, $i =1,\dots,n$ are local coordinates in the base $M$ defined
on $V$, we have local coordinates $(x^i,y^\alpha)$, $i =1,\dots,n$,
$\alpha =1,\dots,\ell$ in $E$. The anchor map and the Lie bracket are
then determined by the local functions $\rho^i_\alpha$ and
$C^\alpha_{\beta\gamma}$ on $M$ (called the \emph{structure functions
  of the Lie algebroid}) defined by
\[
\rho(e_\alpha)= \sum_{i=1}^n \rho^i_\alpha\pd{}{x^i} \qquad \text{and}
\qquad \br{e_\alpha}{e_\beta} = \sum_{\gamma=1}^\ell
C^\gamma_{\alpha\beta}e_\gamma.
\]
The structure functions satisfy the following relations
\begin{align}\label{eq:structure-equations}
  \rho^j_\alpha\pd{\rho^i_\beta}{x^j} -
  \rho^j_\beta\pd{\rho^i_\alpha}{x^j} = \rho^i_\gamma
  C^\gamma_{\alpha\beta} \, , \quad \text{and} \quad
  \sum_{\mathrm{cyclic}(\alpha,\beta,\gamma)} \left[\rho^i_\alpha\pd{
      C^\nu_{\beta\gamma}}{x^i} + C^\mu_{\alpha\nu}
    C^\nu_{\beta\gamma}\right]=0 \, ,
\end{align}
where the summation over repeated indexes is understood.
Equations~\eqref{eq:structure-equations} are usually called the
\emph{structure equations of the Lie algebroid}. Finally, the Lie
bracket of two sections of $E$ can be expressed in terms of the basis
$\{e_\alpha\}_{\alpha=1}^\ell$ as
\begin{align} \label{eq:Lie-bracket}
  \br{\sigma}{\eta}=\left( \sigma^\gamma \rho_\gamma^k
    \pd{\eta^\alpha}{x^k} - \eta^\gamma \rho_\gamma^k
    \pd{\sigma^\alpha}{x^k} + C^\alpha_{\beta\gamma} \sigma^\beta
    \eta^\gamma \right) e_\alpha \, .
\end{align}
If $\famY$ is a family of sections of $E$, we will denote by
$\invclos\famY$ the distribution obtained by closing (the distribution
defined by) $\famY$ under the Lie bracket.

\subsection{Admissible maps and morphisms of Lie
  algebroids}\label{se:morphisms-general}

Let $\map{\tau}{E}{M}$ and $\map{\ov{\tau}}{\ov{E}}{\ov{M}}$ be two
Lie algebroids with associated anchor maps $\map{\rho}{E}{TM}$ and
$\map{\ov{\rho}}{\ov{E}}{T\ov{M}}$.  A bundle map
$\map{\Psi}{E}{\ov{E}}$ is said to be \emph{admissible} if
$T\psi\circ\rho = \ov{\rho}\circ\Psi$. Equivalently $\Psi$ is
admissible if and only if it maps admissible curves into admissible
curves. Indeed, if $a(t)$ is admissible on $E$ and projects to $m(t)$,
then $\ov{a}(t)=\Psi(a(t))$ projects to $\ov{m}(t)=\psi(m(t))$ and it
is admissible, since
\[  
\ov{\rho}(\ov{a}(t)) = \ov{\rho}(\Psi(a(t)) =
T\psi(\rho(a(t))=T\psi(\dot{m}(t)) = \dot{\ov{m}}(t).
\]
Denoting by $\map{\psi}{M}{\ov{M}}$ the map on the base, one has the
following commutative diagram {\let\ds\displaystyle
\[
\xymatrix{%
  &&TM\ar[rr]^{\ds
    T\psi}\ar[lddd]^{\ds\tau_{M}}&&T\ov{M}\ar[lddd]^{\ds\tau_{\ov{M}}}
  \\
  &E\ar[dd]_{\ds\tau}\ar[ru]^{\ds\rho}\ar[rr]^{\ds\Psi} &&
  \ov{E}\ar[dd]_{\ds\ov{\tau}}\ar[ru]^{\ds\ov{\rho}}
  \\
  &&&\\
  &M\ar[rr]_{\ds\psi}&&\ov{M}}
\]

A map $\map{\Psi}{E}{\ov{E}}$ is a \emph{morphism of Lie algebroids}
if it is admissible and preserves the Lie algebra structure of the
algebroids~\cite{PJH-KM:90}, that is, for any $\sigma$ and $\eta$ sections
of $E$ such that there exist some sections $\{\zeta_l\}_{l=1}^p$ of
$\ov{E}$ and some functions $F_l$, $G_l$, $l=1,\dots,p$ on $M$ with
\begin{align*}
  \Psi\circ\sigma=\sum_{l=1}^p F_l (\zeta_l\circ\psi) \, , \quad
  \Psi\circ\eta=\sum_{l=1}^p G_l (\zeta_l\circ\psi) \, ,
\end{align*}
then, the image of the Lie bracket of $\sigma$ and $\eta$ under $\Psi$
is
\[
\Psi\circ[\sigma,\eta]= \sum_{l=1}^p \left(\rho(\sigma)G_l -
  \rho(\eta)F_l \right)\, (\zeta_l\circ\psi) +\sum_{l_1,l_2=1}^p
F_{l_1} G_{l_2} \,([\zeta_{l_1},\zeta_{l_2}]\circ\psi).
\]
Given local basis $\{e_\alpha\}_{\alpha=1}^\ell$ and
$\{\ov{e}_\alpha\}_{\alpha=1}^{\ov{\ell}}$ of sections of $E$ and
$\ov{E}$, respectively, a bundle map $\Psi$ can be written $\Psi
(e_\alpha) = \Psi_\alpha^\beta \ov{e}_\beta$ for certain local
functions $\Psi_\alpha$ on $M$. Then, one can check that $\Psi$ is a
morphism if and only if
\begin{equation}\label{eq:morphism-coordinates}
  \Psi^\beta_\gamma C^\gamma_{\alpha\delta} =
    \left(\rho^i_\alpha\pd{\Psi^\beta_\delta}{x^i} -  
    \rho^i_\delta\pd{\Psi^\beta_\alpha}{x^i}\right) +
  \ov C^\beta_{\theta\sigma}\Psi^\theta_\alpha\Psi^\sigma_\delta \, .
\end{equation}

Notice that if $\sigma$, $\eta$ are $\Psi$-related to sections
$\ov{\sigma}$, $\ov{\eta} \in \Sec{\ov{E}}$, i.e., $\Psi\circ\sigma =
\ov{\sigma} \circ \psi$ and $ \Psi\circ\eta = \ov{\eta} \circ \psi$,
then the Lie bracket $[\sigma,\eta]$ is $\Psi$-related to the Lie
bracket $[\ov{\sigma},\ov{\eta}]$, $\Psi\circ[\sigma,\eta] =
[\ov{\sigma},\ov{\eta}]\circ\psi$.

\subsection{Linear connections}\label{se:linear-connections}

Here we briefly present the notion of $E$-connection on a vector
bundle (cf.~\cite{RLF:02}, see also~\cite{FC-BL:02,AG-EM:03}), and
discuss some related objects.

\begin{definition}
  Let $\tau: E \rightarrow M$ be a Lie algebroid.  A linear
  $E$-connection on a vector bundle $\map{\pi}{P}{M}$ is a
  $\real$-bilinear map $\map{\nabla}{\sec{E}\times\sec{P}}{\sec{P}}$
  such that
  \[
  \nabla_{F\sigma}\alpha = F\nabla_\sigma\alpha \qquad \text{and}
  \qquad \nabla_{\sigma}(F\alpha) = (\rho(\sigma)F)\alpha
  +F\nabla_\sigma\alpha
  \]
  for any function $F\in\cinfty{M}$, section $\sigma$ of $E$ and
  section $\alpha$ of $P$.
\end{definition}

Throughout the paper, we will restrict our attention to the case
$P=E$, and by a \emph{connection on $E$} we will understand a linear
$E$-connection on $\tau:E \rightarrow M$.  Given a local basis
$\{e_\alpha\}_{\alpha=1}^\ell$ of sections of $E$, the local
expression of the covariant derivative is
\[
\nabla_\sigma\eta = \sigma^\alpha \Bigl(
\rho^i_\alpha\pd{\eta^\gamma}{x^i} + \CC^\gamma_{\alpha\beta}
\eta^\beta\Bigr) e_\gamma.
\]
The terms $\CC^\gamma_{\alpha\beta}$ are called the \emph{connection
  coefficients}.  As in the study of tangent bundle geometry, the
skew-symmetric part of the connection defines the so-called
\emph{torsion tensor},
\[
T(\sigma,\eta) = \nabla_\sigma\eta - \nabla_\eta\sigma - [\sigma,\eta]
\, ,
\]
and the symmetric part of the connection determines what we call the
\emph{symmetric product},
\[
\symprod{\sigma}{\eta} = \nabla_\sigma\eta + \nabla_\eta\sigma \, .
\]
The local expression of the symmetric product is (compare with the
expression for the Lie bracket~\eqref{eq:Lie-bracket})
\begin{align}\label{eq:symmetric-product}
  \symprod{\sigma}{\eta} = \left( \sigma^\gamma \rho_\gamma^k
    \pd{\eta^\alpha}{x^k} + \eta^\gamma \rho_\gamma^k
    \pd{\sigma^\alpha}{x^k} + \SC^\alpha_{\beta\gamma} \sigma^\beta
    \eta^\gamma \right) e_\alpha \, , \qquad \text{where} \quad
  \SC^\gamma_{\alpha\beta} = \CC^\gamma_{\alpha\beta}
  +\CC^\gamma_{\beta\alpha} \, .
\end{align}
In particular, notice that
$\symprod{e_\alpha}{e_\beta}=\SC^\gamma_{\alpha\beta} e_\gamma$.
Similarly as with the involutive closure, if $\famY$ is a family of
sections of $E$, we will denote by $\symclos\famY$ the distribution
obtained by closing (the distribution defined by) $\famY$ under the
symmetric product.

Since the covariant derivative is $\cinfty{M}$-linear in the first
argument, it is possible to define the derivative of a section $\sigma
\in \Sec{E}$ with respect to an element $a\in E_m$ by simply putting
\[
\nabla_a\sigma = \nabla_\eta\sigma(m) \, ,
\]
where $\eta \in \Sec{E}$ is any section such that $\eta(m)=a$.
Moreover, the covariant derivative allows us to take the derivative of
sections along maps and, as a particular case, of sections along
curves. If we have a morphism of Lie algebroids $\map{\Phi}{F}{E}$
over the map $\map{\varphi}{N}{M}$, then we can define the derivative
of a section of $E$ along $\varphi$ as follows.

\begin{definition}
  Let $\sigma : N \rightarrow E$ be a section of $E$ along $\varphi$,
  i.e., $\sigma (n) \in E_{\varphi(n)}$, $n \in N$. Then $\sigma$ can
  be written in the form $\sigma = \sum_{l=1}^p F_l
  (\zeta_l\circ\varphi)$, for some sections $\{ \zeta_1, \dots,
  \zeta_p \} \subset \sec{E}$ and some functions $\{F_1,\dots,F_p \}
  \subset \cinfty{N}$.  The \emph{derivative of $\sigma$ along
    $\varphi$} is defined by
  \[
  \nabla_b \sigma = \sum_{l=1}^p \left[ (\rho(b)F_l)
    \zeta_l(\varphi(n))+F_l(n)\nabla_{\Phi(b)} \, \zeta_l \right] \, ,
  \quad b \in F \, .
  \]
\end{definition}

\begin{remark}
  {\rm Within this framework, one can consider time-dependent sections
    of $E$ as follows: take the morphism $\map{\Phi}{T\real\times
      E}{E}$, $\Phi(t,a)=a$ over the map $\map{\varphi}{\real\times
      M}{M}$, $\varphi(t,m)=m$.  The Lie algebroid structure on
    $T\real\times E$ is the direct product structure, that is, the
    anchor is $\rho_{T\real\times E}(\tau\frac{d}{dt},a)=
    \tau\pd{}{t}+\rho(a)$, and the bracket of projectable sections (on
    both factors) is the sum of the brackets on $T\real$ and $E$.  }
\end{remark}

When studying mechanical control systems related by a morphism of Lie
algebroids, we will resort to the following notion concerning the
interplay between maps and linear connections.

\begin{definition}
  Let $\nabla$ and $\ov{\nabla}$ be connections on $E$ and $\ov{E}$,
  respectively, and let $\Psi$ be a bundle map from $E$ to $\ov{E}$.
  We say that $\Psi$ \emph{maps the connection $\nabla$ to the
    connection $\ov{\nabla}$} if
  \[
  \Psi\circ(\nabla_\sigma\eta) = \ov{\nabla}_{\sigma}(\Psi\circ\eta).
  \]
\end{definition}
In coordinates this condition is equivalent to
\begin{equation} \label{eq:mapcon-coordinates}
  \Psi^\beta_\gamma\CC^\gamma_{\alpha\delta} = 
  \rho^i_\alpha\pd{\Psi^\beta_\delta}{x^i}+ 
  \ov\CC^\beta_{\theta\sigma}\Psi^\theta_\alpha\Psi^\sigma_\delta \, .
\end{equation}

\subsubsection*{Geodesics}

Consider the following situation: let $a: t \mapsto E$ be an
admissible curve, and let $b: t \mapsto E$ be a curve in $E$, both of
them projecting by $\tau$ onto the same base curve in $M$, $\tau
(a(t)) = m(t) = \tau (b(t))$. Take the Lie algebroid structure $T\real
\to \real$ and consider the morphism $\Phi : T \real \rightarrow E$,
$\Phi (t,\dot{t})= \dot{t} a(t)$ over $\varphi: \real \rightarrow M$,
$\varphi(t)=m(t)$. Then one can define the \emph{derivative of $b(t)$
  along $a(t)$} as $\nabla_{\frac{d}{dt}}b(t)$.  In the literature,
this derivative is usually denoted by $\nabla_{a(t)}b(t)$.  In local
coordinates, this reads
\[
\nabla_{a(t)}b(t) = \left[
  \frac{db^\gamma}{dt}+\CC^\gamma_{\alpha\beta}a^\alpha b^\beta
\right]e_\gamma(m(t)).
\]

\begin{definition}
  Let $\tau: E \rightarrow M$ be a Lie algebroid and $\nabla$ a
  connection on $E$. An admissible curve $a: t \mapsto E$ is said to
  be a \emph{geodesic of $\nabla$} if $\nabla_{a(t)}a(t)=0$.
\end{definition}
In local coordinates, the conditions for being a geodesic reads
\begin{align}\label{eq:geodesic-equation}
  \frac{da^\gamma}{dt}+\frac{1}{2}\SC^\gamma_{\alpha\beta}a^\alpha
  a^\beta=0 \, .
\end{align}

\subsubsection*{The Levi-Civita
  connection}\label{se:Levi-Civita-connection}

Let $\map{\Gc}{E\times_M E}{\real}$ be a bundle metric on a Lie
algebroid $\tau: E \rightarrow M$. In a parallel way to the situation
in tangent bundle geometry, one can see that there is a canonical
connection on $E$ associated with $\Gc$. The proof is analogous and
will be omitted.

\begin{proposition}
  Given a bundle metric $\Gc$ on $E$, there is a unique connection
  $\nabla^\Gc$ on $E$ which is torsion-less and metric with respect to
  $\Gc$. The connection $\nabla^\Gc$ is determined by the formula
  \begin{align*}
    2\Gc(\nabla_\sigma\eta,\zeta) = \rho(\sigma)\Gc(\eta,\zeta) +
    \rho(\eta)\Gc(\sigma,\zeta) - \rho(\zeta)\Gc(\eta,\sigma) +
    \Gc(\sigma,[\zeta,\eta]) + \Gc(\eta,[\zeta,\sigma])
    -\Gc(\zeta,[\eta,\sigma]) \, ,
  \end{align*}
  for $\sigma, \eta,\zeta \in \Sec{E}$.
\end{proposition}

Denoting by $\{e_\alpha\}_{\alpha=1}^\ell$ a local basis of sections
of $E$, and by $\{e^\alpha\}_{\alpha=1}^\ell$ its dual basis, the
bundle metric can be locally written as $\Gc = \Gc_{\alpha \beta} \,
e^\alpha \otimes e^\beta$. The connection coefficients of
$\nabla^{\Gc}$ are
\[
\CC^\alpha_{\beta\gamma} = \frac{1}{2}\Gc^{\alpha\nu}
\Bigl([\nu,\beta;\gamma] + [\nu,\Gamma;\beta]+[\beta,\gamma;\nu]\Bigr)
\, ,
\]
where $(\Gc^{\mu \nu})$ is the inverse matrix of $(\Gc_{\alpha
  \beta})$, and $[\alpha,\beta;\gamma]$ is a shorthand notation for
\[
[\alpha,\beta;\gamma] = \pd{\Gc_{\alpha\beta}}{x^i}\rho^i_\gamma +
C^\mu_{\alpha\beta}\Gc_{\mu\gamma}.
\]

Associated with the bundle metric $\Gc$, one has the \emph{musical
  isomorphisms}
\begin{align*}
  \flat_\Gc : E \rightarrow E^* \, , \; \pai{\flat_\Gc (a)}{b} = \Gc
  (a,b) \, , \qquad \sharp_\Gc : E^* \rightarrow E \, , \; \sharp_\Gc
  (\theta) = \flat_\Gc^{-1} (\theta) \, ,
\end{align*}
where $E^* \rightarrow M$ denotes the dual bundle of $E$.  Given a
function $V$ on $M$, the \emph{gradient of $V$}, $\grad_{\Gc} V$, is
the section of $E$ defined by $\grad_{\Gc} V = \sharp_{\Gc} (\rho^* dV
)$.

\subsubsection*{The constrained
  connection}\label{se:constrained-connection}

Here, we introduce the notion of constrained connection on a Lie
algebroid, which generalizes the concept of nonholonomic
connection~\cite{ADL:96a,JLS:28}. This notion will be later useful to
model control systems subject to nonholonomic constraints, following
the developments in~\cite{AMB-PEC:98,ADL:97a}. Given an arbitrary
connection and an arbitrary subbundle of a Lie algebroid, one can
define a new connection which enjoys special properties with respect
to the subbundle.  This connection is determined by the choice of a
projector. Let $D$ be a subbundle of $E$ and let $P$ be a projector
onto $D$, $P: E \rightarrow D$. We denote by $Q$ the complementary
projector of $P$, $Q=I-P$, and by $D^c$ the complementary subbundle
$D^c=\textrm{Im}(Q)$. In this way, one has $D\oplus D^c=E$.

\begin{definition}\label{dfn:constrained-connection}
  Given a connection $\nabla$ on $E$ and a projector map $P: E
  \rightarrow D$, the \emph{constrained connection} is the connection
  $\cnabla$ on $E$ defined by
  \begin{align*}
    \cnabla_\sigma\eta=P(\nabla_\sigma\eta)+\nabla_\sigma(Q\eta) \, ,
    \quad \sigma,\eta \in \Sec{E} \, .
  \end{align*}
\end{definition}

Some interesting properties of the constrained connection $\cnabla$
are the following. Their proof is straightforward and will be omitted
for brevity.

\begin{proposition}\label{prop:constrained-connection-properties}
  The following properties of the constrained connection $\cnabla$
  hold:
  \begin{enumerate}
  \item The connection $\cnabla$ restricts to $D$, i.e.,
    $\cnabla_\sigma\eta \in D$ for any $\eta\in\sec{D}$ and
    $\sigma\in\sec{E}$.
  \item The symmetric product $\symprod{\cdot}{\cdot}\check{}$
    associated with $\cnabla$ is given by
    \begin{itemize}
    \item $\symprod{\sigma}{\eta}\check{} = P(\symprod{\sigma}{\eta})$
      for $\sigma, \eta \in \sec{D}$.
    \item $\symprod{\sigma}{\eta}\check{} = P(\symprod{\sigma}{\eta})+
      \symprod{\sigma}{\eta}$ for $\sigma,\eta\in\sec{D^c}$.
    \item $\symprod{\sigma}{\eta}\check{} = P(\symprod{\sigma}{\eta})+
      \nabla_{\sigma}{\eta}$ for $\sigma\in\sec{D}$ and
      $\eta\in\sec{D^c}$.
    \end{itemize}
  \item The torsion tensor $\check{T}$ of $\cnabla$ is given by
    \begin{itemize}
    \item $\check{T}(\sigma,\eta) = P(T(\sigma,\eta)) -
      Q([\sigma,\eta]) $ for $\sigma,\eta\in\sec{D}$.
    \item $\check{T}(\sigma,\eta) = P(T(\sigma,\eta)) +
      \nabla_{\sigma}{\eta} - \nabla_{\eta}{\sigma} -
      Q([\sigma,\eta])$ for $\sigma,\eta\in\sec{D^c}$.
    \item $\check{T}(\sigma,\eta) = P(T(\sigma,\eta)) +
      \nabla_{\sigma}{\eta} - Q([\sigma,\eta])$ for $\sigma\in\sec{D}$
      and $\eta\in\sec{D^c}$.
    \end{itemize}
  \end{enumerate}
\end{proposition}

\begin{proposition}
  Let $\tau: E \rightarrow M$, $\tau : \ov{E} \rightarrow \ov{M}$ be
  two Lie algebroids, with projectors $P: E \rightarrow D$ and $P:
  \ov{E} \rightarrow \ov{D}$.  Let $\nabla$ and $\ov{\nabla}$ be
  connections on $E$ and $\ov{E}$, respectively. Assume that a
  morphism of Lie algebroids $\map{\Psi}{E}{\ov{E}}$ maps $\nabla$
  onto $\ov{\nabla}$. If $\Psi\circ P=\ov{P}\circ\Psi$ (equivalently
  $\Psi(D)\subset \ov{D}$ and $\Psi(D^c)\subset \ov{D}^c$), then
  $\Psi$ maps $\cnabla$ onto $\check{\ov{\nabla}}$.
\end{proposition}
\begin{proof}
  Since $Q=I-P$, it follows that $\Psi\circ Q=\ov{Q}\circ\Psi$.
  Therefore, for all $\eta \in \Sec{E}$ and $b \in E$
  \begin{align*}
    \Psi(\cnabla_b\eta) & = \Psi(P(\nabla_b\eta)) +
    \Psi(\nabla_b(Q\eta)) =
    \ov{P}(\Psi(\nabla_b\eta)) + \ov{\nabla}_b(\Psi\circ (Q\eta))\\
    & =
    \ov{P}(\ov{\nabla}_b(\Psi\circ\eta))+\nabla_b(\ov{Q}(\Psi\circ\eta))
    = \check{\ov{\nabla}}_b(\Psi\circ\eta) \, .
  \end{align*}
\end{proof}

\section{The prolongation of a Lie
  algebroid}\label{se:prolongation} \setcounter{equation}{0}

Here we briefly review the notion of the prolongation of a Lie
algebroid. For further details, see~\cite{EM:01a}. Given a Lie
algebroid $E$, the underlying motivation behind the introduction of
the prolongation of $E$ is that of formulating second-order dynamical
systems on $E$. Thinking of $E$ as a substitute of the tangent bundle
of $M$, the tangent bundle of $E$ is not the appropriate space to
describe second-order dynamics on $E$. This is clear if we note that
the projection to $M$ of a vector tangent to $E$ is a vector tangent
to $M$, and what one would like instead is an element of $E$, the
\emph{new} tangent bundle of $M$.

A space which takes into account this restriction is the $E$-tangent
bundle of $E$, also called the \emph{prolongation of $E$}, which we
denote by $\prol{E}$. This Lie algebroid is defined as the vector
bundle $\map{\tau_1}{\prol{E}}{E}$ whose fiber at a point $a\in E_m$
is the vector space
\[
\prol[a]{E} =\set{(b,v)\in E_m\times T_aE}{\rho(b)=T_a\tau(v)} \,.
\]
Note that if the fibers of $E$ are $\ell$-dimensional, then the fibers of
$\prol{E}$ are $2 \ell$-dimensional.  We will use the redundant notation
$(a,b,v)$ to denote the element $(b,v)\in\prol[a]{E}$.

The anchor of $\prol{E}$ is the map $\map{\rho^1}{\prol{E}}{TE}$,
defined by $\rho^1(a,b,v)=v$. We also consider the map
$\map{\prb}{\prol{E}}{E}$ defined by $\prb(a,b,v)=b$. The Lie bracket
associated with $\prol{E}$ is defined as follows in terms of
projectable sections.  A section $Z$ of $\prol{E}$ is
\emph{projectable} if there exists a section $\sigma$ of $E$ such that
$\prb\circ Z=\sigma\circ\tau$.  Equivalently, a section $Z$ is
projectable if and only if it is of the form $Z(a)=(a,\sigma(\tau
(a)),X(a))$, for some section $\sigma$ of $E$ and some vector field
$X$ on $E$. The Lie bracket of two projectable sections $Z_1$ and
$Z_2$ is then given by
\[
[Z_1,Z_2](a)=(a,[\sigma_1,\sigma_2](m),[X_1,X_2](a)) \, , \quad a \in E \, .
\]
It is easy to see that $[Z_1,Z_2](a)$ is an element of $\prol{E}$ for
every $a\in E$. Since any section of $\prol{E}$ can be locally written
as a linear combination of projectable sections, the definition of the
Lie bracket for sections of $\prol{E}$ follows.

Given local coordinates $(x^i,y^\alpha)$ associated with a basis
$\{e_\alpha\}$ of sections of $E$, we can define a local basis
$\{\X_\alpha,\V_\alpha\}_{\alpha=1}^\ell$ of sections of $\prol{E}$ by
\begin{align}\label{eq:local-basis-prol}
  \X_\alpha(a) =
  \Bigl(a,e_\alpha(\tau(a)),\rho^i_\alpha\pd{}{x^i}\at{a}\Bigr)
  \qquad\text{and} \qquad \V_\alpha(a) =
  \Bigl(a,0,\pd{}{y^\alpha}\at{a}\Bigr) \, , \quad \alpha =1, \dots,
  \ell \, .
\end{align}
If $(a,b,v)$ is an element of $\prol{E}$, with $b=z^\alpha e_\alpha$
and $v=\rho^i_\alpha z^\alpha\pd{}{x^i}+v^\alpha\pd{}{y^\alpha}$, then
we can write
\[
(a,b,v)=z^\alpha\X_\alpha(a)+v^\alpha\V_\alpha(a).
\]
The Lie brackets of the elements of the basis are
\[
\br{\X_\alpha}{\X_\beta}= C^\gamma_{\alpha\beta}\:\X_\gamma \, ,
\qquad \br{\X_\alpha}{\V_\beta}=0 \, ,\qquad
\br{\V_\alpha}{\V_\beta}=0 \, .
\]
Finally, notice that the anchor map $\rho^1$ applied to a section $Z =
Z^\alpha\X_\alpha+V^\alpha\V_\alpha$ of $\prol{E}$ defines a vector
field on $E$ whose coordinate expression is
\[
\rho^1(Z) = \rho^i_\alpha Z^\alpha \pd{}{x^i} + V^\alpha\pd{}{y^\alpha} \, .
\]

The following notion shows that morphisms of Lie algebroids can also be
prolonged.

\begin{definition}
  Let $\Psi :E \rightarrow \ov{E}$ be a an admissible map between two
  Lie algebroids $\tau:E \rightarrow M$ and $\tau:\ov{E} \rightarrow
  \ov{M}$. The \emph{prolongation of $\Psi$} is the mapping
  $\map{\prolmap{\Psi}}{\prol{E}}{\prol{\ov{E}}}$ defined by $
  \prolmap{\Psi}(a,b,v)=(\Psi(a),\Psi(b),T_a\Psi(v))$.
\end{definition}

It is not difficult to see that if $\Psi$ is a morphism of Lie
algebroids, then its prolongation is also a morphism of Lie
algebroids. Given local basis $\{e_\alpha\}_{\alpha=1}^\ell$ and
$\{\ov{e}_\alpha\}_{\alpha=1}^{\ov{\ell}}$ of sections of $E$ and
$\ov{E}$, respectively, if
$\Psi(e_\alpha)=\Psi_\alpha^\beta\ov{e}_\beta$, then the action of the
prolongation of $\Psi$ on $\prol{E}$ is determined by
\begin{equation}  \label{eq:prolongation-coordinates}
  \begin{aligned}
    \prolmap{\Psi}(\X_\alpha(a)) &= \Psi^\beta_\alpha(m)\ov\X_\beta
    (\Psi(a)) +\rho^i_\alpha(m) \pd{\Psi^\gamma_\beta}{x^i}(m)a^\beta
    \ov\V_\gamma(\Psi(a)) \, , \; a \in E_m \, , \\
    \prolmap{\Psi}(\V_\alpha(a)) &= \Psi^\beta_\alpha(m) \, \ov
    \V_\beta (\Psi(a)) \, , \; a \in E_m \, .
  \end{aligned}
\end{equation}

\subsection{Vertical and horizontal
    sections}\label{se:vertical-horizontal} 

An element $(a,b,v) \in \prol{E}$ is said to be \emph{vertical} if it
is of the form $(a,0,v)$, with $v$ a vertical vector tangent to $E$ at
$a$. It follows that the vertical space of $\prol{E}$ at the point
$a\in E_m$, which we denote by $\Ver_{a} (\prol{E})$, can be
identified with $E_m$ by (a slight modification of) the usual vertical
lifting map:
\begin{align*}
  b\in E_m \longmapsto (a,0,b\spV_a) \in \Ver_{a} (\prol{E}) \, ,
\end{align*}
where $b\spV_a\in T_aE$ is the tangent vector to the curve $t\mapsto a
+ t b$ at $t=0$.  If $\sigma$ is a section of $E$, then the section
$\sigma\spV$ of $\prol{E}$ defined by
$\sigma\spV(a)=(a,0,\sigma(m)\spV_a)$ will be called the
\emph{vertical lift of~$\sigma$}. Vertical elements are linear
combinations of $\{\V_\alpha\}_{\alpha=1}^\ell$.  Specifically, if the
section $\sigma$ of $E$ has the local expression $\sigma=\sigma^\alpha
e_\alpha$, then $\sigma\spV$ is of the form
\[
\sigma\spV=\sigma^\alpha\V_\alpha.
\]

Consider the \emph{zero-section of $\tau:E \lo M$}, that is
\[
0_M : M \lo E \, , \quad 0_M(m) = 0_m \, .
\]
This section is a canonical embedding of $M$ into $E$. Consequently,
we can regard $TM$ as a subspace of $TE$. Now, define the
\emph{horizontal space} along $0_M$,
\[
\Hor_m (\prol{E}) = \set{ (0_m,b,v) \in \prol[0_m]{E}}{v \in T_mM \subseteq
  T_{0_m} E} , \qquad m \in M \, .
\]
Note that $\tau_1 (\Hor (\prol{E})) = \mathrm{Im}(0_M)$, i.e.
$\Hor(\prol{E})$ is only defined along points of $E$ that belong to the
zero-section $0_M$. Moreover, $\Hor_m(\prol{E})$ is a vector subspace of
$\prol[0_m]{E}$. The following result shows that the horizontal space is
complementary to the vertical space along the zero-section $0_M$.

\begin{lemma}\label{le:decomposition}
  Along the zero-section of $\tau:E \rightarrow M$, we have the following
  direct sum decomposition
  \[
  \prol[0_M]{E} = \Hor(\prol{E}) \oplus \Ver_{0_M} (\prol{E}) \, .
  \]
\end{lemma}

\begin{proof}
  Consider the map $\map{\Hc}{E}{\prol[0_M]{E}}$ given by
  $\Hc(b)=(0_{\tau (b)},b,T0_M(\rho(b)))$. The image of $\Hc$ is
  precisely the horizontal space, $\Hor_m(\prol{E}) = \Hc(E_m)$, $m
  \in M$. Note that $\Hc$ is a splitting of the short exact sequence
  \[
  \xymatrix{0\ar[r]&E\ar[r]^{V} & \prol[0_M]{E} \ar[r]^{\prb} &
    E\ar[r] & 0} ,
  \]
  where $V (b) = (0_{\tau (b)},0,b\spV_{0_{\tau(b)}})$. As a
  consequence, we have $\prol[0_m]{E} = \Hc (E_m) \oplus \Ver_{0_m}
  (\prol{E})$, $m \in M$.  The restriction of $\prb$ to
  $\Hor(\prol{E})$ is an isomorphism $\map{\prb}{\Hor(\prol{E})}{E}$,
  whose inverse map is $\Hc$.  For $(0_m,b,v) \in \prol[0_m]{E}$, the
  decomposition is given by $(0_m,b,v) = (0_m,b, T (0_M \circ \tau)
  (v)) + (0_m,0_m,v-T (0_M \circ \tau) (v)) \in \Hor(\prol{E}) +
  \Ver_{0_M} (\prol{E})$.
\end{proof}

\subsection{Homogeneity}\label{se:homogeneity}

A property that plays an important role in our later analysis is that
of homogeneity. Consider the section $\Delta$ of $\prol{E}$ defined by
$\Delta(a)=(a,0,a\spV_a)$, $a \in E$. The section $\Delta$ is called
the \emph{Liouville section of $\prol{E}$}. In coordinates, we have
\[
\Delta = y^\alpha\V_\alpha \qquad \text{and} \qquad
\rho^1(\Delta)=y^\alpha\pd{}{y^\alpha}.
\]

A function $F \in \cinfty{E}$ is said to be \emph{homogeneous of
  degree $s\in \integer$} if
\[
\mathcal{L}_{\rho^1 (\Delta)} F = s F \, .
\]
where $\mathcal{L}$ stands for the Lie derivative
operator\footnote{Alternatively, one can see that homogeneous
  functions of degree $s$ verify $d_\Delta F=s F$. (See~\cite{PL-CMM:87}
  for the precise definition of the derivative operator $d_\Delta$).}.
In a local chart, a homogeneous function of degree $s\geq0$ is a
homogeneous polynomial in $\{y^\alpha \}_{\alpha=1}^\ell$ of degree
$s$ with arbitrary functions of $(x^i)_{i=1}^n$ as coefficients.
Consequently, homogeneous functions of degree $0$ are (pullbacks of)
functions on the base $M$ and there are not (smooth) non-trivial
functions homogeneous of degree $s\leq-1$. A section $Z$ of $\prol{E}$
is said to be \emph{homogeneous of degree $s\in \integer$} if
\[
[\Delta,Z]=s Z.
\]
We denote by $\Pc_s$ the set of homogeneous sections of $\prol{E}$ of
degree $s$. The following result describes the basic properties
concerning homogeneous sections. The proof is omitted for brevity.

\begin{lemma}\label{le:facts-homogeneity}
  Let $r$, $s \in \integer$ and let $Z$ be a section of $\prol{E}$.
  Then
  \begin{itemize}
  \item[(i)] $[\Pc_s,\Pc_r] \subseteq \Pc_{s+r}$, and $\Pc_{s}=\{0\}$
    if $s\leq2$,
  \item[(ii)] $Z\in\Pc_{-1}$ if and only if there exists a section
    $\sigma$ of $E$ such that $Z = \sigma\spV$,
  \item[(iii)] $Z\in\Pc_{0}$ if and only if $Z$ is a projectable
    section,
  \item[(iv)] $\X_\alpha \in \Pc_{0}$ and $\V_\alpha \in \Pc_{-1}$,
    for $\alpha =1, \dots, \ell$. Moreover, if $Z = Z^\alpha \X_\alpha
    + V^\alpha \V_\alpha$ is the local expression of $Z$, then $Z \in
    \Pc_{s}$ if and only if the functions $Z^\alpha$ are homogeneous
    of degree $s$ and the functions $V^\alpha$ are homogeneous of
    degree $s+1$.
  \end{itemize}
\end{lemma}

Note that for all $Z\in \Pc_s$, $s\ge 1$, we have $Z (0_m)=0_{0_m}$,
$m \in M$, that is, the homogeneous sections of $\prol{E}$ of degree
greater or equal than $1$ vanish at the zero-section of $E$.

\subsection{SODE sections}\label{se:sodes} 

The dynamics of a mechanical system evolving on a certain
configuration manifold is described by means of a vector field on the
tangent bundle of the manifold which is a second order differential
equation.  Likewise, we will need the notion of \sode\ section of the
prolongation to describe the behavior of mechanical systems evolving
on Lie algebroids. In this section, we introduce this concept and
discuss several geometric properties, resembling those of second order
differential equations on manifolds.

\begin{definition}
  A section $\Gamma$ of $\prol{E}$ is a second-order differential
  equation (\sode) section on the Lie algebroid $E$ if
  $\prol{\tau}\circ\Gamma=\textrm{Id}_E$.
\end{definition}

A vector $v \in T_a E$ is called \emph{admissible} if $T \rho (v) =
\rho(a)$. Note that a curve in $E$ is admissible if and only if its
tangent vectors are admissible. We denote by $\Adm(E)$ the set of all
admissible tangent vectors. Notice that $v$ is admissible if and only
if $(a, a, v) \in \prol{E}$.  Therefore we can identify $\Adm (E)$
with the subset of $\prol{E}$ formed by all elements of that form,
\begin{align*}
  \Adm (E) = \setdef{z \in \prol{E}}{\tau_1 (z) = \prol{\tau}(z)} \, .
\end{align*}
Equivalently, a \sode\ section can be defined as a section of
$\prol{E}$ that takes values in $\Adm (E)$, i.e.,
$\Gamma(a)=(a,a,X(a))$, $a\in E$, where $X$ is a vector field on $E$
verifying $\rho(a)=T_a\tau(X(a))$. If $\Gamma$ is a \sode\ section,
then it has the local expression $\Gamma = y^\alpha\X_\alpha +
F^\alpha(x,y)\V_\alpha$ and its associated vector field is of the form
$$
\rho^1(\Gamma) = \rho^i_\alpha
y^\alpha\pd{}{x^i}+F^\alpha(x,y)\pd{}{y^\alpha}.
$$
The integral curves of this vector field satisfy the differential
equations
\[
\dot{x}^i = \rho^i_\alpha y^\alpha \, , \qquad \dot{y}^\alpha =
F^\alpha(x,y) \, .
\]
In particular, we are specially interested in homogeneous \sode\ 
sections with degree $1$. We will refer to such \sode\ sections as
\emph{sprays}. Locally, a spray is such that the functions $F^\alpha$
are homogeneous with degree 2, $F^\alpha(x,y) =
-\frac{1}{2}\SC^\alpha_{\beta\gamma}(x)y^\beta y^\gamma$ for some
symmetric coefficients $\SC^\alpha_{\beta\gamma}$.  The sprays are in
one to one correspondence with torsion-less $E$-covariant derivatives
(cf. Section~\ref{se:linear-connections}), as we show in the
following.

Let $\nabla$ be a connection on $E$.  The geodesic
equations~\eqref{eq:geodesic-equation}, together with the
admissibility condition, correspond to the differential equations of
the integral curves of a \sode\ section $\Gamma_{\nabla}$ of
$\prol{E}$, which is locally given by
\[
\Gamma_{\nabla} = y^\alpha\X_\alpha - \frac{1}{2} \left(
  \CC^\alpha_{\beta\gamma} +\CC^\alpha_{\gamma\beta} \right) y^\beta
y^\gamma \V_\alpha \, .
\]
From the coordinate expression, we easily see that $\Gamma_{\nabla}$
is homogeneous with degree one, i.e., it is a spray. Moreover, notice
that $\Gamma_\nabla$ is determined by the symmetric product associated
with $\nabla$ (cf.  eq.~\eqref{eq:symmetric-product}).

We now show how a spray $\Gamma$ and a (2,1)-tensor field $T$
determine a connection on $E$ in a unique way. As an intermediate
step, we first define the symmetric product associated with a spray.
Given $\sigma,\eta\in\sec{E}$, consider the section of $\prol{E}$,
$[\sigma\spV,[\Gamma,\eta\spV]]$.  From the properties of the Lie
bracket, we deduce that this section is (smooth and) homogeneous with
degree $-1$. From Lemma~\ref{le:facts-homogeneity}(ii), it follows
that $[\sigma\spV,[\Gamma,\eta\spV]]$ is the vertical lift of a
section of $E$, which we denote by $\symprod{\sigma}{\eta}_{\Gamma}$.
Therefore,
\[
\symprod{\sigma}{\eta}\spV=[\sigma\spV,[\Gamma,\eta\spV]].
\]
From the Jacobi identity, one can deduce that the operation
$\symprod{\sigma}{\eta}$ is symmetric.  Moreover, if $f$ is a function
on $M$, we have $\symprod{\sigma}{f\eta} =
\rho(\sigma)f\eta+f\symprod{\sigma}{\eta}$, since
\[
\symprod{\sigma}{f\eta}\spV = [[ \sigma\spV,\Gamma],(f\eta)\spV] =
\rho^1([\sigma\spV,\Gamma])f\eta\spV + f[\sigma\spV,[\Gamma,\eta\spV]]
= (\rho(\sigma)f\eta+f\symprod{\sigma}{\eta})\spV \, .
\]

\begin{definition}
  Given a spray $\Gamma$ on $E$, $\symprod{\cdot}{\cdot}_{\Gamma}$ is
  the \emph{symmetric product associated with $\Gamma$}.
\end{definition} 

Taking a local basis of sections of $E$, one can compute
\begin{align}
  [\Gamma,\sigma\spV] & = [y^\alpha\X_\alpha-\frac{1}{2}
  \SC^\alpha_{\beta\gamma}y^\beta
  y^\gamma\V_\alpha ,\sigma^\beta\V_\beta] 
  = -\sigma^\alpha\X_\alpha + y^\alpha\left(
    \rho^i_\alpha\pd{\sigma^\gamma}{x^i}+
    \SC^\gamma_{\alpha\beta}\sigma^\beta
  \right) \V_\gamma, \label{eq:grado-0} \\
  [\eta\spV,[\Gamma,\sigma\spV]] &= \left(
    \sigma^\gamma\rho_\gamma^k\pd{\eta^\alpha}{x^k}
    +\eta^\gamma\rho_\gamma^k\pd{\sigma^\alpha}{x^k} +
    \SC^\alpha_{\beta\gamma} \sigma^\beta \eta^\gamma \right)
  \V_\alpha \, .
  \label{eq:grado-1}
\end{align}
From the local expression of the symmetric product, one can see that
the symmetric product \emph{determines} and \emph{is determined} by
the spray $\Gamma$.

\begin{proposition}
  Given a spray $\Gamma$ and a skew-symmetric (2,1) tensor field $T$
  on $E$, there exists a unique connection
  $\map{\nabla^{\Gamma,T}}{\sec{E}\times\sec{E}}{\sec{E}}$ on $E$ such
  that its associated spray is $\Gamma$ and its torsion is~$T$. This
  connection is given by
  \begin{align}\label{eq:associated-connection}
    \nabla^{\Gamma,T}_\sigma\eta = \frac{1}{2}\big([\sigma,\eta] +
    T(\sigma,\eta)\big) + \frac{1}{2}\symprod{\sigma}{\eta}_\Gamma \,.
  \end{align} 
\end{proposition}
\begin{proof}
  From the properties of the symmetric product and the Lie bracket, it
  follows that $\nabla^{\Gamma,T}$ defined
  in~\eqref{eq:associated-connection} is a connection on $E$. Its
  torsion is given by
  \[
  \nabla^{\Gamma,T}_\sigma\eta -
  \nabla^{\Gamma,T}_\eta\sigma-[\sigma,\eta] =
  \frac{1}{2}\big([\sigma,\eta] + T(\sigma,\eta)\big) -
  \frac{1}{2}\big([\eta,\sigma] + T(\eta,\sigma)\big) - [\sigma,\eta]
  = T(\sigma,\eta).
  \]
  By definition, the spray associated with $\nabla^{\Gamma,T}$ is the
  \sode\ section $Z$ whose projection $\rho^1(Z)$ is determined by the
  differential equation $\nabla^{\Gamma,T}_{a(t)}a(t)=0$, for
  admissible curves $a: \real \rightarrow E$, $t \mapsto a(t)$.  A
  simple coordinate calculation shows that these equations can be
  locally written as
  $\dot{y}^\alpha+\frac{1}{2}\SC^\alpha_{\beta\gamma}y^\beta
  y^\gamma=0$, which concludes the proof.
\end{proof}

The connection coefficients of $\nabla^{\Gamma,T}$ are given by
\begin{align*}
  \CC^\alpha_{\beta\gamma} = \frac{1}{2}\Big(\SC^\alpha_{\beta\gamma}
  + T^\alpha_{\beta\gamma} +C^\alpha_{\beta\gamma}\Big).
\end{align*}
Note also that the symmetric product associated with $\Gamma$
precisely corresponds to the symmetric product defined by
$\nabla_{\Gamma,T}$.
\[
\symprod{\sigma}{\eta}_{\Gamma} = \nabla^{\Gamma,T}_\sigma\eta +
\nabla^{\Gamma,T}_\eta\sigma \, .
\]

\begin{proposition}
  Let $\tau:E \rightarrow M$ and $\tau:\ov{E} \rightarrow \ov{M}$ be
  two Lie algebroids. Let $\nabla$ and $\ov{\nabla}$ be connections on
  $E$ and $\ov{E}$, respectively. A morphism $\Psi$ maps $\nabla$ onto
  $\ov{\nabla}$ if and only if $\Psi$ maps the associated spray
  $\Gamma$ into the associated spray $\ov{\Gamma}$ and maps the
  torsion tensor $T$ into the torsion tensor $\ov{T}$.
\end{proposition}

\begin{proof}
  We prove in coordinates the above statement for torsion-less
  connections, from where the general result follows easily. If we
  take a \sode\ $\Gamma=y^\alpha\X_\alpha+F^\alpha\V_\alpha$, then
  \[
  \prolmap{\Psi}(\Gamma(a))= \ov{a}^\alpha\ov\X_\alpha+\left[
    \rho^i_\alpha(m)a^\alpha\pd{\Psi^\gamma_\beta}{x^i}(m)a^\beta
    +F^\alpha(a)\Psi^\gamma_\alpha(m) \right] \ov\X_\gamma \, , a \in
  E_m \, .
  \]
  This last expression is equal to $\ov\Gamma(\ov{a}) =
  \ov{a}^\alpha\ov\X_\alpha+\ov F^\alpha(\ov{a})\ov\V_\alpha$ if and
  only if
  \[
  \ov F^\alpha(\ov{a})=\Psi^\gamma_\alpha F^\alpha(a)+a^\beta
  \rho^i_\alpha a^\alpha\pd{\Psi^\gamma_\beta}{x^i} \, .
  \]
  In the case of two sprays, $F^\gamma =
  -\frac{1}{2}\SC^\gamma_{\alpha\beta}y^\alpha y^\beta$ and $\ov
  F^\gamma=-\frac{1}{2}\ov\SC^\gamma_{\alpha\beta}\ov y^\alpha \ov
  y^\beta$, so that the above equation reads
  \[
  \Psi^\beta_\gamma\SC^\gamma_{\alpha\delta} =
  \left(\rho^i_\alpha\pd{\Psi^\beta_\delta}{x^i}+
    \rho^i_\delta\pd{\Psi^\beta_\alpha}{x^i}\right)+
  \ov\SC^\beta_{\theta\sigma}\Psi^\theta_\alpha\Psi^\sigma_\delta.
  \]
  Since $\Psi$ is a morphism, it verifies
  equation~\eqref{eq:morphism-coordinates}. Summing both expressions,
  we finally obtain equation~\eqref{eq:mapcon-coordinates}, as
  claimed.
\end{proof}

\subsection{Geodesically invariant
  subbundles}\label{se:geodesically-invariant}

Here we introduce the notion of geodesically invariant subbundles
(which is a generalization of the concept of geodesically invariant
distributions~\cite{ADL:96a}), and establish its relation with the
symmetric product associated with the connection.

\begin{definition}\label{dfn:geodesically-invariant}
  Let $\nabla$ be a  connection on $E$.  A subbundle $D\subset
  E$ is \emph{geodesically invariant for $\nabla$} if every geodesic
  $\map{a}{[t_0,t_1]}{E}$ such that $a(t_0)\in D$ verifies that
  $a(t)\in D$ for all $t\in[t_0,t_1]$.
\end{definition}

\begin{proposition}\label{prop:geodesic-invariance}
  Let $\nabla$ be a  connection on $E$.  A subbundle $D\subset
  E$ is geodesically invariant for $\nabla$ if and only if it is
  invariant under the symmetric product, i.e.
  $\symprod{\sigma}{\eta}\in\sec{D}$ for every
  $\sigma,\eta\in\sec{D}$.
\end{proposition}
\begin{proof}
  By definition, the subbundle $D$ is geodesically invariant if and
  only if the spray of the connection $\Gamma_{\nabla}$ is tangent to
  $D$, i.e. $\rho^1(\Gamma)|_D\in TD$. The latter is equivalent to the
  condition $\left( {\cal L}_{\rho^1 (\Gamma)}\phi \right)_{|D}=0$ for
  every function $\phi$ on $E$ such that $\phi_{|D}=0$.  Since $D$ is
  linear, the constraint functions $\phi$ are of the form
  $\phi=\hat{\mu}$, with $\mu\in\sec{D^\circ}$ and
  $\hat{\mu}(a)=\pai{\mu_m}{a}$, $a\in E_m$.  Therefore, $D$ is
  geodesically invariant if and only if $({\cal L}_{\rho^1
    (\Gamma)}\hat{\mu})\circ\sigma=0$ for all $\mu\in\sec{D^\circ}$
  and $\sigma\in\sec{D}$. Let us see that $({\cal L}_{\rho^1
    (\Gamma)}\hat{\mu})\circ\sigma =
  -\frac{1}{2}\pai{\mu}{\nabla_\sigma\sigma}$.  In coordinates
  \[
  {\cal L}_{\rho^1 (\Gamma)}\hat{\mu}={\cal L}_{\rho^1
    (\Gamma)}(\mu_\alpha y^\alpha)
  =\rho^i_\beta\pd{\mu_\alpha}{x^i}y^\alpha y^\beta
  -\frac{1}{2}\mu_\gamma\SC^\gamma_{\alpha\beta}y^\alpha y^\beta \, .
  \]
  Taking the restriction to the image of $\sigma$, and writing
  explicitly the symmetric parts, we have
  \[
  {\cal L}_{\rho^1 (\Gamma)}\hat{\mu}\circ\sigma = \frac{1}{2}\left(
    \rho^i_\beta\pd{\mu_\alpha}{x^i} +
    \rho^i_\alpha\pd{\mu_\beta}{x^i} -
    \mu_\gamma\SC^\gamma_{\alpha\beta} \right)
  \sigma^\alpha\sigma^\beta \, .
  \]
  On the other hand
  \[
  \nabla_\sigma\mu = \sigma^\alpha\left(
    \rho^i_\alpha\pd{\mu_\beta}{x^i} - \mu_\gamma
    \CC^\gamma_{\beta\alpha}\right)e^\beta \, ,
  \]
  and thus
  \begin{align*}
    \pai{\nabla_\sigma\mu}{\sigma} = \sigma^\alpha \sigma^\beta\left(
      \rho^i_\alpha\pd{\mu_\beta}{x^i} -
      \mu_\gamma\CC^\gamma_{\beta\alpha}\right) =
    \frac{1}{2}\sigma^\alpha \sigma_\beta\left(
      \rho^i_\alpha\pd{\mu_\beta}{x^i} + \rho^i_\beta
      \pd{\mu_\alpha}{x^i} -\mu_\gamma
      \SC^\gamma_{\beta\alpha}\right),
  \end{align*}
  where we have used the fact that $\CC^\gamma_{\alpha\beta} +
  \CC^\gamma_{\beta\alpha} = \SC^\gamma_{\beta\alpha}$.  Therefore
  $({\cal L}_{\rho^1 (\Gamma)}\hat\mu)\circ\sigma = \frac{1}{2}
  \pai{\nabla_\sigma\mu}{\sigma}$, and the result follows by taking
  into account that $\pai{\nabla_\sigma\mu}{\sigma} =
  \rho(\sigma)\pai{\mu}{\sigma} - \pai{\mu}{\nabla_\sigma\sigma} =
  -\pai{\mu}{\nabla_\sigma\sigma}$, because $\pai{\mu}{\sigma}=0$.
\end{proof}

Given a connection $\nabla$ on $E$ and a projector map $P: E
\rightarrow D$, consider the constrained connection $\cnabla$
introduced in Section~\ref{se:constrained-connection}.  Using the
above result and
Proposition~\ref{prop:constrained-connection-properties}(ii), one can
deduce that $D$ is geodesically invariant for $\cnabla$.

\section{General control systems on Lie algebroids}
\label{se:general-control-systems}
\setcounter{equation}{0}

In this section we present the notion of a control system on a Lie
algebroid.  We introduce the concept of accessibility subbundle in the
Lie algebroid and provide basic tests for controllability, building on
the known results for control systems defined on manifolds. Finally,
we study the controllability properties of control systems related by
means of a morphism of Lie algebroids.

Consider a Lie algebroid $\tau:E \lo M$, with anchor map $\rho:E \lo
TM$. Let $\sigma$, $\eta_1,\dots,\eta_k$ be sections of $E$. A
\emph{control problem on the Lie algebroid $E \lo M$ with drift
  section $\sigma$ and input sections $\eta_1,\dots,\eta_k$} is
defined by the following equation on $M$,
\begin{equation}\label{eq:control-system}
  \dot{m} (t) = \rho\Bigl(\sigma (m(t))) + \sum_{i=1}^k u_i(t)\eta_i
  (m(t))\Bigr) \, ,
\end{equation}
where $u=(u_1,\dots, u_k) \in U$, and $U$ is an open set of $\real^k$
containing $0$. The function $ t \mapsto u(t) = (u_1(t), \ldots,
u_m(t))$ belongs to a certain class of functions of time, denoted by
$\mathcal{U}$, called the \emph{set of admissible controls}.  For our
purposes, we may restrict the admissible controls to be the piecewise
constant functions with values in $U$. Notice that the trajectories of
the control system are admissible curves of the Lie algebroid, and
therefore they must lie on a leaf of $E$.  It then follows that if $E$
is not transitive, then there are points that cannot be connected by
solutions of any control system defined on such a Lie algebroid.  In
particular, the system~\eqref{eq:control-system} cannot be locally
accessible at points $m \in M$ where $\rho$ is not surjective.  Since
the emphasis here is put on the controllability analysis, without loss
of generality we will restrict our attention to locally transitive Lie
algebroids.
  
Denoting by $f=\rho(\sigma)$ and $g_i=\rho(\eta_i)$, we can rewrite the
system~\eqref{eq:control-system} as
\begin{equation}\label{eq:standard}
  \dot{m} (t) = f (m(t)) + \sum_{i=1}^k u_i(t) g_i (m(t)) \, ,
\end{equation}
which is a standard nonlinear control system on $M$ affine in the
inputs.  Here we make use of the additional geometric structure
provided by the Lie algebroid in order to carry over the analysis of
the controllability properties of the control
system~\eqref{eq:control-system}. We refer to~\cite{HN-AJvdS:90} for a
comprehensive discussion of the notions of reachable sets,
accessibility algebra and computable accessibility tests. A short list
of definitions is provided in the appendix for reference.
  
\begin{definition}
  The \emph{accessibility algebra~$\D$ of the control
    system~\eqref{eq:control-system} in the Lie algebroid} is the
  smallest subalgebra of $\Sec{E}$ that contains the sections
  $\sigma,\eta_1,\dots,\eta_k$.
\end{definition}

Using the Jacobi identity, one can deduce that any element of
accessibility algebra~$\D$ is a linear combination of repeated Lie
brackets of sections of the form
\[
[\zeta_l,[\zeta_{l-1},[\dots,[\zeta_2,\zeta_1]\dots]]] \, ,
\]
where $\zeta_i \in \{\sigma,\eta_1,\dots,\eta_k\}$, $1 \le i \le l$ and $l
\in \natural$.

\begin{definition}
  The \emph{accessibility subbundle in the Lie algebroid}, denoted by
  $\invclos{\{\sigma,\eta_1,\dots,\eta_k\}}$, is the vector subbundle of $E$
  generated by the accessibility algebra~$\D$,
  \[
  \invclos{\{\sigma,\eta_1,\dots,\eta_k\}} = \spn \left\{ \zeta(m) \;
    | \; \zeta \; \text{section of $E$ in} \; \D \right\} \, , \quad m
  \in M \, .
  \]
\end{definition}
If the dimension of $\invclos{\{\sigma,\eta_1,\dots,\eta_k\}}$ is
constant, then $\invclos{\{\sigma,\eta_1,\dots,\eta_k\}}$ is the
smallest Lie subalgebroid of $E$ that has
$\{\sigma,\eta_1,\ldots,\eta_k\}$ as sections.  In the following
result, we establish accessibility tests for control systems of the
form~\eqref{eq:control-system}.

\begin{theorem}\label{th:accessibility-test-in-E}
  Consider the system~\eqref{eq:control-system}.  Let $m \in M$ and
  assume the Lie algebroid $E$ is locally transitive at $m$. Then,
  $\invclos{\{\sigma,\eta_1,\dots,\eta_k\}}(m) + \ker\rho(m)= E_{m}$
  implies that the system is locally accessible from $m$.
\end{theorem}

\begin{proof}
  We first prove that every element $X$ of the accessibility algebra $\C$
  of~\eqref{eq:standard} can be written as $X=\rho (\zeta)$, with $\zeta \in
  \D$ (see the appendix for the precise definition of the accessibility
  algebra~$\C$).  Take an element of $\C$ of the form
  $[X_l,[X_{l-1},[\dots,[X_2,X_1]\dots]]]$, with $X_i \in \{f = \rho
  (\sigma),g_1=\rho (\eta_1),\dots,g_k = \rho (\eta_k) \}$. Denote $X_i =
  \rho (\zeta_i)$, with $\zeta_i \in \{\sigma,\eta_1,\dots,\eta_k\}$. Since
  $\rho$ is a Lie algebra homomorphism, we have that
  \[
  [X_l,[X_{l-1},[\dots,[X_2,X_1]\dots]]] = \rho \left(
    [\zeta_l,[\zeta_{l-1},[\dots,[\zeta_2,\zeta_1]\dots]]] \right) \, .
  \]
  Since $\rho$ is linear, we conclude the accessibility subbundle in the Lie
  algebroid is mapped by $\rho$ onto the accessibility distribution $C$.
  Now, we show $C(m) = T_{m} M$ if and only if
  $\invclos{\{\sigma,\eta_1,\dots,\eta_k\}}(m) + \ker\rho(m)= E_{m}$.  Assume
  first that $C(m) = T_{m} M$.  Let $e_{m} \in E_{m}$. Consider $\rho (e) \in
  T_{m}M = C(m) = \rho (\invclos{\{\sigma,\eta_1,\dots,\eta_k\}}(m))$.  There
  exists $a \in \invclos{\{\sigma,\eta_1,\dots,\eta_k\}}(m)$ such that $\rho
  (e) = \rho (a)$.  Then, $e-a \in \ker \rho (m)$ and $e = a + (e-a)$.
  Therefore $E_{m} = \invclos{\{\sigma,\eta_1,\dots,\eta_k\}}(m) + \ker \rho
  (m)$.  In addition, $\rho (E_{m}) = \rho
  (\invclos{\{\sigma,\eta_1,\dots,\eta_k\}}(m)) = T_{m}M$. The other
  implication is trivial.  Finally, the result follows from Chow's
  theorem~\cite{HN-AJvdS:90}.
\end{proof}

\begin{remark}
  {\rm If the dimension of $\invclos{\{\sigma,\eta_1,\dots,\eta_k\}}$
    is constant, then the above theorem expresses the following fact:
    if the Lie subalgebroid $\invclos{\{\sigma,\eta_1,\dots,\eta_k\}}$
    is locally transitive at $m$, then the
    system~\eqref{eq:control-system} is locally accessible from $m$.}
\end{remark}

In practice, the most interesting property to establish is
controllability.  In what follows, we provide a controllability test
for systems of the form~\eqref{eq:control-system} by adapting the
notions of good and bad Lie brackets of vector fields~\cite{HJS:87} to
the setting of Lie algebroids.  Let $\{X_0,X_1, \dots, X_k\}$ be a set
of sections of the Lie algebroid $E$.  The \emph{degree} of an
iterated Lie bracket $B$ of elements in $\{X_0,X_1, \dots, X_k\}$ is
the number of occurrences of all its factors, and is therefore given
by $\delta (B) = \delta_0(B) + \delta_1 (B) + \dots + \delta_k(B)$,
where $\delta_i(B)$ is the number of times that $X_i$ appears in $B$.
A Lie bracket $B$ is said to be \emph{bad} if $\delta_0(B)$ is odd and
$\delta_i(B)$ is even, $i \in \{ 1, \dots, k \}$.  Otherwise, $B$ is
said to be \emph{good}. To make precise sense of these notions
(degree, bad, good) one must resort to the concept of free Lie
algebras, but it should be clear from the context what we mean here
(see~\cite{HJS:87} for a detailed discussion).

\begin{theorem}\label{th:controllability-test-in-E}
  Assume that the system~\eqref{eq:control-system} is locally
  accessible from $m \in M$. If every bad Lie bracket $B$ in $\{
  \sigma, \eta_1, \dots, \eta_k \}$ evaluated at $m$ can be put as an
  $\real$-linear combination of good Lie brackets in $\{ \sigma,
  \eta_1, \dots, \eta_k \}$ of lower degree than~$B$ and elements in
  the kernel of the anchor map $\rho$ at $m$, $\ker \rho (m)$, then
  the system is locally controllable from $m$.
\end{theorem}

\begin{proof}
  Under the hypothesis of the theorem, any bad Lie bracket of the
  vector fields $\{ f=\rho(\sigma), g_1=\rho(\eta_1), \dots,
  g_k=\rho(\eta_k) \}$ evaluated at $m$ can be written a
  $\real$-linear combination of good Lie brackets in $\{ f, g_1,
  \dots, g_k \}$ of lower degree. The application of Theorem~7.3
  in~\cite{HJS:87} gives the result.
\end{proof}

\subsection{Control systems related by a morphism of Lie algebroids} 

Let $\tau :E \rightarrow M$ and $\ov{\tau}:\ov{E} \rightarrow \ov{M}$
be two Lie algebroids. Consider a control system defined on $E$ with
drift section $\sigma$ and independent input sections $\eta_i$,
$i=1,\ldots,k$, and a control system defined on $\ov{E}$ with drift
section $\ov{\sigma}$ and independent input sections $\ov{\eta}_j$,
$j=1,\ldots,\ov{k}$.  Let $\Psi : E \rightarrow \ov{E}$ be a morphism
of Lie algebroids.  Let $\spn \{ \ov{\eta}_1,\dots \ov{\eta}_{\ov{k}}
\}$ denote the linear subbundle of $\ov{E}$ generated by the sections
$\ov{\eta}_j$, $j=1,\ldots,\ov{k}$, and consider the affine subbundle
$\xi + \spn \{ \ov{\eta}_1,\dots \ov{\eta}_{\ov{k}} \}$.

\begin{definition}
  The control problems on $E$ and $\ov{E}$ are \emph{weakly
    $\Psi$-related} (or $\Psi$ maps the system on $E$ onto the system
  on $\ov{E}$) if $\Psi (\sigma (m)) \in \xi (\psi(m))+ \spn \{
  \ov{\eta}_1,\dots \ov{\eta}_{\ov{k}} \} (\psi (m))$ and $\Psi
  (\eta_i (m)) \in \spn \{ \ov{\eta}_1,\dots \ov{\eta}_{\ov{k}} \}
  (\psi (m))$ for all $i=1,\dots, k$ and all $m \in M$.
\end{definition}

Equivalently, two systems are weakly $\Psi$-related if there exist
functions $C^i_j$ and $b^j$ on $M$ such that
\begin{align}\label{eq:related-systems-weakly}
  \Psi\circ\sigma = \xi\circ\psi + \sum_{j=1}^{\ov{k}}
  b^j(\ov{\eta}_j\circ\psi) \quad \text{and} \quad
  \Psi\circ\eta_i = \sum_{j=1}^{\ov{k}} C^j_i(\ov{\eta}_j\circ\psi) \,
  , \quad i=1,\dots ,k \, .
\end{align}

\begin{definition}
  The control problems on $E$ and $\ov{E}$ are \emph{$\Psi$-related} if
  $\sigma$ is $\Psi$-related to a section of $\ov{E}$ with values in $\xi +
  \spn \{ \ov{\eta}_1,\dots \ov{\eta}_{\ov{k}} \}$ and $\eta_i$ is
  $\Psi$-related to a section of $\ov{E}$ with values in $ \spn \{
  \ov{\eta}_1,\dots \ov{\eta}_{\ov{k}} \}$, for all $i=1,\dots, k$.
\end{definition}

Equivalently, two systems are $\Psi$-related if there exist functions
$\ov{C}^i_j$ and $\ov{b}^j$ on $\ov{M}$ such that
\begin{align}\label{eq:related-systems}
  \Psi\circ\sigma = \xi\circ\psi + \sum_{j=1}^{\ov{k}} (\ov{b}^j
  \ov{\eta}_j)\circ\psi \quad \text{and} \quad
  \Psi\circ\eta_i =
  \sum_{j=1}^{\ov{k}} (\ov{C}^j_i \ov{\eta}_j )\circ\psi \, , \quad
  i=1,\dots ,k \, .
\end{align}
Clearly, $\Psi$-related systems are also weakly $\Psi$-related. The
following result establishes the relation between the controllability
properties of weakly $\Psi$-related systems.

\begin{proposition}\label{th:open-map-general}
  Let $\map{\Psi}{E}{\ov{E}}$ be a morphism of Lie algebroids such
  that the associated base map $\psi$ is open. Consider two control
  systems on $E$ and $\ov{E}$ that are weakly $\Psi$-related. If the
  system on $E$ is locally accessible (respectively locally
  controllable) from $m \in M$, then the system on $\ov{E}$ is locally
  accessible (respectively locally controllable) from $\psi(m)$.
\end{proposition}

\begin{proof}
  First, note that if two control systems are $\Psi$-related, then the
  image by the morphism $\Psi$ of a solution $a(t)$ of the system on
  $E$ with control functions $t \mapsto u_i(t)$, $i=1,\dots,k$ is a
  solution of the system on $\ov{E}$ with control functions $t \mapsto
  \ov{u}_j(t) = b^j(\tau(a(t)))+\sum_{i=1}^k C^j_i(\tau(a(t)))u_i(t)$,
  $j=1,\dots,\ov{k}$.
  
  Let $m \in M$ and take $\ov{V}$ an open neighborhood of $\psi(m)$ in
  $\ov{M}$.  Consider the open neighborhood $V=\psi^{-1}(\ov{V})$ of
  $m$. If the system on $E$ is locally accessible from $m$, then there
  is a non-empty open set $\calo$ contained in $\R_M^{V}(m,\leq T)$
  (where $\R_M^V(m,\le T)$ denotes the reachable set in $M$ starting
  from $m$ in time less than or equal to $T$, see the appendix).
  Since $\psi$ is an open map and $\psi(\R_M^V(m,t)) \subset
  \R_{\ov{M}}^{\psi(V)}(\psi(m),t)$, it follows that
  $\ov{\calo}=\psi(\calo)\subset\R_{\ov{M}}^{\ov{V}}(\psi(m),\leq T)$,
  and thus the system on $\ov{E}$ is locally accessible (respectively
  controllable) from $\psi(m)$. The argument for the controllable case
  is analogous.
\end{proof}

A interesting particular case occurs when $\Psi$ is an isomorphism
between the fibers of the Lie algebroids.  In such a case, if two
systems are $\Psi$-related, one can see that the sufficient conditions
for local accessibility (respectively controllability) are either
simultaneously satisfied on $E$ and $\ov{E}$ or simultaneously not
satisfied.

\begin{proposition}\label{prop:equivalence-general}
  Let $\Psi:E \rightarrow \ov{E}$ be a morphism of Lie algebroids
  which is an isomorphism on each fiber. Consider two control systems
  on $E$ and $\ov{E}$, with $k \ge \ov{k}$, that are $\Psi$-related.
  Let $m \in M$. Then
  \begin{itemize}
  \item[(a)] $\invclos{\{\sigma,\eta_1,\dots,\eta_k\}}(m) = E_m$ if
    and only if
    $\invclos{\{\ov{\sigma},\ov{\eta}_1,\dots,\ov{\eta}_{\ov{k}}\}}
    (\psi (m)) = \ov{E}_{\psi (m)}$,
  \item[(b)] Every bad Lie bracket in $\{\sigma,\eta_1,\dots,\eta_k\}$
    evaluated at $m$ can be put as an $\real$-linear combination of
    good Lie brackets of lower degree if and only if every bad Lie
    bracket in $\{\ov{\sigma},\ov{\eta}_1,\dots,\ov{\eta}_{\ov{k}}\}$
    evaluated at $\psi(m) \in \ov{M}$ can be put as an $\real$-linear
    combination of good Lie brackets of lower degree.
  \end{itemize}
\end{proposition}

\begin{proof}  
  Since the control systems are $\Psi$-related and $\Psi$ is an
  isomorphism on each fiber, we have that $k = \ov{k}$, and one can
  assume without loss of generality that the drift sections $\sigma$
  and $\ov{\sigma}$ and the input sections $\eta_i$ and $\ov{\eta}_i$
  satisfy
  \begin{equation}\label{eq:related-by-psi-general}
    \Psi \circ \sigma = \ov{\sigma} \circ \psi 
    \qquad\text{and}\qquad
    \Psi \circ \eta_i = \ov{\eta}_i \circ \psi \, , \quad i=1,\ldots,k \, .
  \end{equation}
  From the properties of a morphism of Lie algebroids (see
  Section~\ref{se:morphisms-general}), we deduce that the Lie bracket
  of any subset of sections in $\{\sigma,\eta_1,\ldots,\eta_k\}$ is
  $\Psi$-related to the Lie bracket of the corresponding subset of
  sections in $\{\ov{\sigma},\ov{\eta}_1,\ldots,\ov{\eta}_k\}$. Using
  now the fact that both Lie algebroids have fibers of the same
  dimension, we conclude the result.
\end{proof}

As an example application, consider the case of a control system on a
manifold $M$ invariant under the action of a symmetry Lie group $G$ on
$M$.  Assume that the action of $G$ on $M$ is free and proper, so that
$\map{\pi}{M}{M/G=\ov{M}}$ is a principal fiber bundle.  Then the
quotient map $\map{\Psi}{TM}{TM/G}$ is a morphism of Lie algebroids
between $E=TM \rightarrow M$ and $\ov{E}=TM/G \rightarrow T(M/G)$.
Moreover, $\Psi$ is an isomorphism in every fiber and its associated
base map $\psi$ is open.  Being the control system on $E$ invariant
under the action of $G$, it induces a control system on $\ov{E}$. From
the above results, we conclude that if the reduced system satisfies
the sufficient conditions for local controllability (respectively
accessibility), then the original system also satisfies such
conditions.  Moreover, since the map $\psi$ is open, if the reduced
system is not locally controllable (respectively accessible), then the
original system cannot be locally controllable (respectively
accessible).

\section{Mechanical control systems}
\label{se:mechanical-control-systems}
\setcounter{equation}{0}

In this section, we consider control problems defined on the
prolongation $\prol{E}$ of $E$ and we make use of the special geometry
of this Lie algebroid to further investigate the controllability
properties of the control systems defined on it.  Given a Lagrangian
function $L: E \rightarrow \real$ on the Lie algebroid, define the
associated action functional ${\cal J} = \int_{t_0}^{t_1}L\,dt$.
Consider the following constrained variational problem: find the
extremals of ${\cal J}$ among the set of admissible curves with fixed
endpoints $m_0$ and $m_1$ in the base $M$, i.e.  curves $a:[t_0,t_1]
\rightarrow E$, $m(t)=\tau(a(t))$, satisfying
\[
\rho(a(t))=\dot{m}(t) \, , \quad m(t_0)=m_0 \quad \text{and} \quad
m(t_1) = m_1 \, .
\]
One can see that the infinitesimal variations $W: [t_0,t_1]
\rightarrow \prol{E}$ of a curve $a:[t_0,t_1] \rightarrow E$
corresponding to this constrained variational problem are of the form
\[
W (t) = \rho^i_\alpha(m(t))\sigma^\alpha(t)\pd{}{x^i} +
\left(\frac{d\sigma^\alpha}{dt} +
  C^\alpha_{\beta\gamma}a^\beta(t)\sigma^\gamma(t)\right)\pd{}{y^\alpha},
\]
where $\sigma(t)$ is a curve on $E$ over $m(t)$ which vanishes at
$t_0$ and $t_1$. A simple calculation and the application of the
Fundamental Theorem of Calculus show that the equations of motion
describing the solutions of the constrained variational problem are
\begin{align*}
  \frac{d}{dt}\pd{L}{y^\alpha} +
  C^\gamma_{\alpha\beta}y^\beta\pd{L}{y^\gamma} =
  \rho^i_\alpha\pd{L}{x^i}.
\end{align*}

\paragraph{Euler-Lagrange operator.}

In alternative terms~\cite{JFC-EM:01,EM:01b,AW:96}, the infinitesimal
variations of a curve $a:[t_0,t_1] \rightarrow E$ can be written as
the restriction to the curve of the complete lift of a general
time-dependent section of $E$. In such a case, the above equations are
precisely the components of the Euler-Lagrange operator $\delta L :
\Adm (E) \rightarrow E^*$, which locally reads
\[
\delta L = \left(\frac{d}{dt}\pd{L}{y^\alpha} +
  C^\gamma_{\alpha\beta}y^\beta\pd{L}{y^\gamma} -
  \rho^i_\alpha\pd{L}{x^i}\right)e^\alpha,
\]
where $\{e^\alpha\}$ is the dual basis of $\{e_\alpha\}$.  The
equations of motion just read
\begin{align*}
  \delta L=0 \, .
\end{align*}
Equivalently, an admissible curve $a:[t_0,t_1] \rightarrow E$ is an
extremal of the action functional ${\cal J}$ if the Euler-Lagrange
operator $\delta L$ vanishes at the points of the curve in $\Adm (E)$,
$t \mapsto (a(t),a(t),\dot{a}(t))$.

\paragraph{Nonholonomic constraints.}

In the case of a system subject to (linear) nonholonomic constraints,
in addition to the above data there is also a subbundle $D$ of the Lie
algebroid which prescribes the allowed velocities for the system.  The
equations of motion then are given by the application of the
Lagrange-d'Alembert principle,
\[
\delta L\in D^\circ \, , \qquad a\in D.
\]
If there is a projector $P: E \rightarrow D$ onto $D$, denoting
$Q=I-P$, the above equations can be rewritten as
\begin{align*}
  P^*(\delta L)=0 \, , \qquad Q(a)=0 .
\end{align*}
Here $P^*$ stands for the dual linear map $\pai{P^*(\theta)}{a} =
\pai{\theta}{P(a)}$, $a \in E$, $\theta \in E^*$.

\paragraph{Control forces.}

In the presence of external forces, the equations of motion for both
the unconstrained and the constrained situations have to be modified.
Assume that some input forces $\{ \theta_1,\dots,\theta_m\} \subset
\sec{E^*}$ act on the Lagrangian system on $E$. Then, the equations of
motion for the unconstrained control system read
\begin{align}\label{eq:Euler-Lagrange-control}
  \delta L = \sum_{l=1}^m u_l \, \theta_l \, ,
\end{align}
and the equations of motion for the nonholonomically constrained
control system are
\begin{align}\label{eq:Lagrange-Alembert-control}
  P^*(\delta L) = \sum_{l=1}^m u_l P^*(\theta_l) \, , \quad Q(a)=0 \, .
\end{align}

\paragraph{Mechanical control systems.}

For the remainder of the paper, we will focus our attention on the
class of Lagrangian control systems $(L,\{\theta_1,\dots,\theta_k\})$
whose Lagrangian function $L : E \rightarrow \real$ is of the form
\[
L(a)=\frac{1}{2} \Gc(a,a) - V \circ \tau (a) \, , \quad a\in E,
\]
with $\map{\Gc}{E\times_M E}{\real}$ a bundle metric on $E$ and $V$ a
function on $M$. We denote by $\{\eta_1,\dots,\eta_k\}$ the input
sections of $E$ determined by the control forces
$\{\theta_1,\dots,\theta_k\}$ via the metric, i.e., $\eta_i =
\sharp_{\Gc} (\theta_i)$.  If $\Gamma_{\nabla^{\Gc}}$ denotes the
spray associated with the Levi-Civita connection $\nabla^\Gc$, the
controlled equations~\eqref{eq:Euler-Lagrange-control} can be written
as
\begin{equation}\label{eq:LCS}
  \dot{a}(t) = \rho^1 \Bigl (\Gamma_{\nabla^\Gc} (a(t)) - (\grad_\Gc
  V) \spV (a(t))+ \sum_{i=1}^k u_i(t) \eta_i\spV (a(t))\Bigr). 
\end{equation}
Note that this system is a control problem on the Lie algebroid
$\prol{E} \rightarrow E$ as defined in
Section~\ref{se:general-control-systems}. This is the reason why we
will refer to the control problem with data
$(\Gc,V,\{\theta_1,\dots,\theta_k\})$ as a \emph{mechanical control
  system on a Lie algebroid}.  Locally, the equations can be written
as
\begin{align*}
  \dot{x}^i & =\rho^i_\alpha y^\alpha , \\[-7pt]
  \dot{y}^\alpha &= -\frac{1}{2} \left( \CC^\alpha_{\beta \gamma} (x)
    +\CC^\alpha_{\gamma\beta} (x) \right) y^\beta y^\gamma -
  \Gc^{\alpha \beta} \rho^i_{\beta} \pd{V}{x^i} +\sum_{i=1}^k
  u_i(t)\eta^\alpha_i(x).
\end{align*}
Alternatively, one can describe the dynamical behavior of the
mechanical control system by means of an equation on $E$ via the
covariant derivative.  An admissible curve $a: t \mapsto a(t)$ is a
solution of the system~\eqref{eq:LCS} if and only if
\begin{align}\label{eq:eqs-motion-connection}
  \nabla^{\Gc}_{a(t)}a(t) + \grad_\Gc V(m(t)) = \sum_{i=1}^k u_i(t)
  \eta_i(m(t)) \, .
\end{align}

\paragraph{Mechanical control systems with constraints.}

If the mechanical control system $(\Gc,V,\{\theta_1,\dots,\theta_k\})$
is subject to the constraints determined by a subbundle $D$ of $E$, we
can do the following.  Consider the orthogonal decomposition
$E=D\oplus D^\perp$, and the associated orthogonal projectors $P:E
\rightarrow D$, $Q:E \rightarrow D^\perp$. Using the fact that
$\Gc(P\cdot,\cdot) = \Gc(\cdot,P\cdot)$, one can write the controlled
equations~\eqref{eq:Lagrange-Alembert-control} as
\[
P(\nabla^\Gc_{a(t)}a(t)) + P(\grad_\Gc V(m(t))) = \sum_{i=1}^k u_i (t)
P(\eta_i(m(t))) \, , \qquad Q(a)=0 \, .
\]
In terms of the constrained connection $\cnabla_\sigma\eta =
P(\nabla^\Gc_\sigma\eta)+\nabla^\Gc_\sigma(Q\eta)$ (cf.
Section~\ref{se:constrained-connection}), we can rewrite this equation
as $\cnabla_{a(t)}a(t) + P(\grad_\Gc V(m(t))) = \sum_{i=1}^k u_i (t)
P(\eta_i(m(t)))$, $Q(a)=0$.  Since the subbundle $D$ is geodesically
invariant for the connection $\cnabla$, it follows that any integral
curve of the spray $\Gamma_{\cnabla}$ associated with $\cnabla$
starting from $a_0 \in D$ is entirely contained in $D$ (cf.
Definition~\ref{dfn:geodesically-invariant}). Since the forcing terms
in~\eqref{eq:eqs-motion-connection-nh} coming from the potential and
the inputs belong to $D$, the same property holds for the total
controlled dynamics. As a consequence, the controlled equations can be
simply stated as
\begin{align}\label{eq:eqs-motion-connection-nh}
  \cnabla_{a(t)}a(t) + P(\grad_\Gc V(m(t))) = \sum_{i=1}^k u_i (t)
  P(\eta_i(m(t))) \, ,\qquad a_0 \in D \, .
\end{align}
Note that one can write the controlled dynamics as a control system on
the Lie algebroid $\prol{E} \rightarrow E$,
\begin{equation}\label{eq:LACS}
  \dot{a}(t) = \rho^1 \Bigl (\Gamma_{\cnabla} (a(t)) - P(\grad_\Gc
  V) \spV (a(t))+ \sum_{i=1}^k u_i(t) P(\eta_i)\spV (a(t))\Bigr). 
\end{equation}
The coordinate expression of these equations is greatly simplified if
we take a basis $\{e_\alpha\}=\{e_a,e_A\}$ of $E$ adapted to the
orthogonal decomposition $E=D\oplus D^\perp$, i.e., $D = \spn
\{e_a\}$, $\D^\perp = \spn \{e_A\}$. Denoting by
$(y^\alpha)=(y^a,y^A)$ the induced coordinates, the constraint
equations $Q(a)=0$ just read $y^A=0$. The controlled
equations~\eqref{eq:eqs-motion-connection-nh} are then
\begin{align*}
  \dot{x}^i & = \rho^i_a y^a , \\
  \dot{y}^a & = - \frac{1}{2}\SC^a_{bc}y^by^c - \Gc^{a\beta}
  \rho^i_\beta\pd{V}{x^i} + \sum_{i=1}^k u_i (t) P(\eta_i)^a , \\
  y^A &= 0.
\end{align*}

\paragraph{Connection control systems.}

Given that the structure of the controlled equations is the same both
in the absence~\eqref{eq:eqs-motion-connection} and in the presence of
constraints~\eqref{eq:eqs-motion-connection-nh}, we will in general
talk about \emph{connection control systems on $\tau:E \rightarrow
  M$}. The dynamics of these systems is governed by an equation of the
type
\begin{align}\label{eq:connection-control}
  \nabla_{a(t)}a(t) + \eta (m(t))) = \sum_{i=1}^k u_i (t) \eta_i
  (m(t))) \, .
\end{align}
Here $\nabla$ is a connection on $E$, and
$\{\eta,\eta_1,\dots,\eta_k\}$ are sections of $E$. We will often
refer to $\eta$ as the potential energy term in
equations~\eqref{eq:connection-control}. Associated with this
equation, there is always a control system on the Lie algebroid
$\prol{E} \rightarrow E$ given by
\begin{equation}\label{eq:CS}
  \dot{a}(t) = \rho^1 \Bigl ( (\Gamma_{\nabla} - \eta\spV )(a(t)) +
  \sum_{i=1}^k u_i(t) \eta_i\spV (a(t))\Bigr).  
\end{equation}

\subsection{Accessibility and controllability
  notions}\label{se:notions}

Here we introduce the notions of accessibility and controllability
that are specialized to mechanical control systems on Lie algebroids.
Let $m \in M$ and consider a neighborhood $V$ of $m$ in $M$. Define
the set of reachable points in the base manifold $M$ starting from $m$
as
\begin{multline*}
  \R_M^V (m, T) = \left\{ m' \in M \, \right. | \left. \exists u \in
    \mathcal{U} \; \hbox{defined on $[0,T]$ such that the evolution
      of~\eqref{eq:CS}} \right. \\
  \left. \hbox{for $a(0)=0_{m}$ satisfies} \; \tau(a(t)) \in V , \, t
    \in [0,T] \, \hbox{and} \; \tau(a(T))=m' \right\} .
\end{multline*}
Alternatively, one may write $\R_M^V (m, T) = \tau
(\R_E^{\tau^{-1}(V)} (0_m, T))$. Denote
\[
\R_M^V (m,\le T) = \bigcup_{t \le T} \R_M^V (m, t) \, .
\]

\begin{definition}\label{dfn:base-controllable}
  The system~\eqref{eq:CS} is \emph{locally base accessible from $m$}
  (respectively, \emph{locally base controllable from $m$}) if
  $\R_M^V(m,\le T)$ contains a non-empty open set of $M$
  (respectively, $\R_M^V(m,\le T)$ contains a non-empty open set of
  $M$ to which $m$ belongs) for all neighborhoods $V$ of $m$ and all
  $T>0$. If this holds for any $m \in M$, then the system is called
  \emph{locally base accessible} (respectively, \emph{locally base
    controllable}).
\end{definition}

In addition to the notions of base accessibility and base
controllability, we shall also consider full-state accessibility and
controllability starting from points of the form $0_m \in E$, $m \in
M$ (note that full-state is meant here with regards to $E$, not to
$TM$).

\begin{definition}\label{dfn:zero-controllable}
  The system~\eqref{eq:CS} is \emph{locally accessible from $m$ at
    zero} (respectively, \emph{locally controllable from $m$ at zero})
  if $\R_E^W(0_{m},\le T)$ contains a non-empty open set of $E$
  (respectively, $\R_E^W(0_{m},\le T)$ contains a non-empty open set
  of $E$ to which $0_{m}$ belongs) for all neighborhoods $W$ of
  $0_{m}$ in $E$ and all $T>0$. If this holds for any $m \in M$, then
  the system is called \emph{locally accessible at zero}
  (respectively, \emph{locally controllable at zero}).
\end{definition}

The relevance of the above definitions stems from the fact that,
frequently, one needs to control a system by starting at rest.
Nevertheless it is important to notice that not every equilibrium
point at $m$ corresponds to the point $0_m$. Indeed, there might be
other relative equilibrium points, explicitly all those points $a \in
E$ such that $\rho^1 (\Gamma (a)) = 0$, i.e., $a$ is in the kernel of
the anchor map $\rho$ and $F^\alpha(a_m)=0$, $\alpha =1,\dots,\ell$.

Finally, we also introduce the notion of accessibility and
controllability with regards to a manifold.

\begin{definition}\label{dfn:fiber-controllable}
  Let $\psi:M \rightarrow N$ be an open mapping.  The
  system~\eqref{eq:CS} is \emph{locally base accessible from $m$ with
    regards to $N$} (respectively, \emph{locally base controllable
    from $m$} with regards to $N$) if $\psi(\R_M^V(m,\le T))$ contains
  a non-empty open set of $N$ (respectively, $\psi(\R_M^V(m,\le T))$
  contains a non-empty open set of $N$ to which $\psi(m)$ belongs) for
  all neighborhoods $V$ of $m$ and all $T>0$.  If this holds for any
  $m \in M$, then the system is called \emph{locally base accessible
    with regards to $N$} (respectively, \emph{locally base
    controllable with regards to $N$}).
\end{definition}

Note that base accessibility and controllability with regards to $M$
with $\Id{M} : M \rightarrow M$ corresponds to the notions of base
accessibility and controllability (cf.
Definition~\ref{dfn:base-controllable}). Moreover, if the system is
base accessible, then it is base accessible with regards to $N$. The
analogous implication for base controllability also holds true.

\subsection{The structure of the control Lie
  algebra}\label{se:structure}

The aim of this section is to show that the analysis of the structure
of the control Lie algebra of affine connection control systems
carried out in~\cite{ADL-RMM:95c} can be further extended to control
systems defined on a Lie algebroid. The enabling technical notion
exploited here is that of homogeneity. As we will show later, this
analysis will allow us to enlarge the class of systems to which the
accessibility and controllability tests can be applied.

For the purpose of evaluating the brackets of the accessibility
subbundle $\invclos{\{\Gamma - \eta\spV, \eta\spV_1, \dots,
  \eta\spV_k\}}$ at initial states of the form $0_m$, $m \in M$, the
discussion on the geometry of $\prol{E}$ along the zero-section (cf.
Section~\ref{se:vertical-horizontal}) will be most helpful. Since the
Lie brackets in the accessibility subbundle of the mechanical control
system are linear combinations of the brackets of the elements
$\{\Gamma, \eta\spV_1, \dots, \eta\spV_k, \eta_0\spV\}$, as an
intermediate step we will first analyze the structure of the subbundle
$\invclos{\{\Gamma, \eta\spV_1, \dots, \eta\spV_k, \eta\spV\}}$.

Let $B$ be a Lie bracket formed with sections of the family $\famX =
\{ \Gamma,\eta_1\spV,\dots,\eta_k\spV,\eta\spV \}$. For each $l$,
consider the following sets
\begin{align*}
  \Br^l (\famX) & = \setdef{B \; \text{bracket in} \, \famX}{\delta
    (B) = l} \, , \\
  \Br_l (\famX) & = \setdef{B \; \text{bracket in} \, \famX}{B \in
    \Pc_l} \, .
\end{align*}

The notion of \emph{primitive} bracket will also be useful. Given a
bracket $B$ in $\famX$, it is clear that we can write $B=[B_1,B_2]$,
with $B_i$ brackets in $\famX$. In turn, we can also write $B_\alpha =
[B_{\alpha1},B_{\alpha2}]$ for $\alpha=1,2$, and continue these
decompositions until we end up with elements belonging to $\famX$.
The collection of brackets $B_1,B_2,B_{11},B_{12},\dots$ are called
the \emph{components} of $B$. The components of $B$ which do not admit
further decompositions are called \emph{irreducible}. A bracket $B$ is
called \emph{primitive} if all of its components are brackets in
$\Br_{-1}(\famX) \cup \Br_{0}(\famX) \cup \{ \Gamma \}$.

We may now recall the following lemma taken from~\cite{ADL:95a}. Although
there it is stated for vector fields, the proof can be readily
extended to sections of a Lie algebroid since it only relies on two
facts: (i) the Jacobi identity, and (ii) the fact that $\famX$ has the
property that any bracket in $\Br_l(\famX)$, $l \le 2$, is identically
zero (which follows from Lemma~\ref{le:facts-homogeneity}(i)).

\begin{lemma}
  Any bracket in $\Br_0(\famX) \cup \Br_{-1}(\famX)$ is a finite sum
  of primitive brackets.
\end{lemma}

As a consequence of this lemma, and the fact that all the brackets in
$\Br_l (\famX)$, with $l\ge 1$ vanish when evaluated at the
zero-section of $E$, we conclude that the only brackets we need to
consider are the primitive brackets in $\Br_{-1}(\famX) \cup
\Br_{0}(\famX)$. We do this next. First, observe the computation of
the basic brackets~\eqref{eq:grado-0} and~\eqref{eq:grado-1}.  In
particular, notice that $[\Gamma,\sigma\spV]$ projects to $-\sigma$.
Second, from Lemma~\ref{le:facts-homogeneity}(ii), we deduce that any
bracket in $\Br_{-1}(\famX)$ is the vertical lift of a section of $E$.
From Lemma~\ref{le:facts-homogeneity}(iii), we deduce that the
brackets $B$ belonging to $\Br_0(\famX)$ are projectable sections.  We
will denote by $\sigma_B$ the section to which it projects. Thus, we
have $B(0_m)=\sigma_B (m)$, $m \in M$. The following result completely
unveils the structure of these brackets.

\begin{lemma}\label{le:primitive-0}
  Let $B\in \Br_{0}(\famX)$ be a primitive bracket. Then either one of
  the following is true,
  \begin{itemize}
  \item[(i)] $B=[\Gamma,B_1]$ with $B_1 \in \Br_{-1}(\famX)$. If
    $\sigma_1$ is the section of $E$ such that $B_1=\sigma_1\spV$,
    then $\sigma_B=-\sigma_1$.  In addition, $[\sigma_2\spV,B] =
    \symprod{\sigma_1}{\sigma_2}\spV$ for all $\sigma_2 \in \Sec{E}$.
  \item[(ii)] $B=[B_1,B_2]$ with $B_1,B_2 \in \Br_0(\famX)$. Then,
    $B(0_m) = [\sigma_{B_1},\sigma_{B_2}](m)$ for all $m \in M$.
  \end{itemize}
\end{lemma}

\begin{proof}
  Let $B \in \Br_{0}(\famX)$. Then, either $B=[\Gamma,B_1]$ with $B_1
  \in \Br_{-1}(\famX)$ primitive, or $B=[B_1,B_2]$ with $B_1,B_2 \in
  \Br_0(\famX)$ both primitive.  In the first case,
  eq.~\eqref{eq:grado-0} gives
  $B(0_m)=[\Gamma,B_1](0_m)=-\sigma_1(m)$, where $B_1=\sigma_1\spV$.
  From~\eqref{eq:grado-1}, we also deduce that
  $[\sigma_2\spV,B]=[\sigma_2\spV,[\Gamma,\sigma_1\spV]] =
  \symprod{\sigma_1}{\sigma_2}\spV$.  In the second case, we have that
  $B_1$ and $B_2$ are projectable onto $\sigma_{B_1}$ and
  $\sigma_{B_2}$ respectively, and therefore $[B_1,B_2]$ projects to
  $[\sigma_{B_1},\sigma_{B_2}]$. Consequently $B(0_m) = [B_1,B_2](0_m)
  = [\sigma_{B_1},\sigma_{B_2}](m)$.
\end{proof}

\begin{proposition}
  \label{prop:decomposition-accessibility-distribution-intermediate} 
  Let $m\in M$. Then,
  \begin{align*}
    & \invclos{\{\Gamma,\eta\spV_1,\dots,\eta\spV_k, \eta\spV\}}
    \cap \Ver_{0_m}(\prol{E}) = \symclos{\{
      \eta,\eta_1,\dots,\eta_k\}}(m)\spV \, , \\
    & \invclos{\{\Gamma,\eta\spV_1,\dots,\eta\spV_k, \eta\spV \}}
    \cap \Hor_{m}(\prol{E}) = \invclos{\symclos{\{ \eta,
        \eta_1,\dots,\eta_k\}}}(m) \, .
  \end{align*}
\end{proposition}

\begin{proof}
  We prove the inclusion $\supseteq$ in the first equality by
  induction. Let us denote $\symcloss{(1)}{\{ \eta,
    \eta_1,\dots,\eta_k\}}=\spn \{ \eta, \eta_1,\dots,\eta_k\}$ and
  \[
  \symcloss{(l)}{\{ \eta, \eta_1,\dots,\eta_k\}} = \spn \{
  \symprod{\sigma_1}{\sigma_2} \; | \; \sigma_i \in
  \symcloss{(l_i)}{\{ \eta,\eta_1,\dots,\eta_k\}} \, , \; l_1 + l_2
  = l \} \, .
  \]
  The result is trivially true for $l=1$, $\symcloss{(1)}{\{ \eta,
    \eta_1,\dots,\eta_k\}}\spV \subseteq
  \invclos{\{\Gamma,\eta\spV_1,\dots,\eta\spV_k,\eta\spV\}} $.
  Assume it is true for $l$ and let us prove it for $l+1$. Take
  $\sigma_i \in \symcloss{(l_i)}{\{ \eta,\eta_1,\dots,\eta_k\}}$,
  $i=1,2$, with $l_1 + l_2 = l+1$ and observe that,
  using~\eqref{eq:grado-1},
  \[
  \symprod{\sigma_1}{\sigma_2}\spV =
  [[\sigma_1\spV,\Gamma],\sigma_2\spV] \, .
  \]
  By induction hypothesis, $\sigma_i\spV \in
  \invclos{\{\Gamma,\eta\spV_1,\dots,\eta\spV_k,\eta\spV\}}$, and
  therefore, $\symprod{\sigma_1}{\sigma_2}\spV \in
  \invclos{\{\Gamma,\eta\spV_1,\dots,\eta\spV_k,\eta\spV\}}$.  From
  the previous discussion, we know that for the opposite inclusion it
  is sufficient to look at the primitive brackets in
  $\Br_{-1}(\famX)$.  Let $B = [B_1,B_2] \in \Br_{-1}(\famX)$
  primitive, with $B_1 \in \Br_{-1}(\famX)$, $B_2 \in \Br_{0}(\famX)$
  primitive. From Lemma~\ref{le:primitive-0}, we have either $B_2 =
  [\Gamma,B_2']$ or $B_2=[B_2',B_2'']$. In the first case, we have
  \[
  B=[B_1,[\Gamma,B_2']] \, .
  \]
  Using again~\eqref{eq:grado-1}, we conclude that $B \in
  \symclos{\{\eta, \eta_1,\dots,\eta_k\}}\spV$.  As for the second
  case, using the Jacobi identity, we have
  \[
  B = [B_1,[B_2',B_2'']] = - [B_2'',[B_1,B_2']] + [B_2',[B_1,B_2'']]
  \, .
  \]
  Applying repeatedly the above argument to $[B_1,B_2']$ and
  $[B_1,B_2'']$ until they are expressed in terms of symmetric
  products, we see that $B$ can be expressed as a linear combination
  of elements in $\symclos{\{ \eta,\eta_1,\dots,\eta_k\}}\spV$, and
  hence we conclude the result.  The second equality is a direct
  consequence of the first one and Lemma~\ref{le:primitive-0}.
\end{proof}

Now, consider the set $\famX' = \{\Gamma - \eta\spV_0,
\eta\spV_1,\dots,\eta\spV_k \}$. As noted before, the elements in
$\invclos{\famX'}$ are linear combinations of the elements in
$\invclos{\famX}$. In fact, for each bracket $B'$ of elements in
$\famX'$, let us define the subset $S(B') \subset \Br(\famX)$ formed
by all possible brackets $B \in \Br (\famX)$ obtained by replacing
each occurrence of $\Gamma-\eta\spV_0$ in $B'$ by either $\Gamma$ or
$\eta\spV_0$. Then, one can prove by induction (cf.~\cite{ADL:95a})
that
\begin{align}\label{eq:decomposition}
  B' = \sum_{B \in S(B')} (-1)^{\delta_{k+1} (B)} B \, ,
\end{align}
where recall that $\delta_{k+1} (B)$ stands for the number of
occurrences of $\eta\spV$ in $B$. Reciprocally, given an element $B
\in \Br (\famX)$, one can determine the bracket $B'$ of elements in
$\famX'$ such that $B \in S(B')$ simply by substituting each
occurrence of $\Gamma$ or $\eta\spV_0$ in $B$ by $\Gamma-\eta\spV_0$.
We denote this operation by $\psinv (B) = B'$.

For each $k \in \natural$, define the following families of sections
in $E$,
\begin{align*}
  \mathcal{C}^{(k)}_{\ver} (\eta;\eta_1,\dots,\eta_k) & =
  \setdef{\sigma \in \Sec{E}}{\sigma\spV = B'', B'' = \hspace*{-.25cm}
    \sum_{
      \begin{subarray}{c}
        \tilde{B} \in S(\psinv(B)) \\
        \cap\Br_{-1}(\famX) \cap \Br_{0}(\famX)
      \end{subarray}
    } \hspace*{-.75cm} (-1)^{\delta_{k+1}(\tilde{B})} \tilde{B} \, ,
    \;  B \in \Br^{2k-1}(\famX) \; \text{primitive}} \, ,
  \\
  \mathcal{C}^{(k)}_{\hor} (\eta;\eta_1,\dots,\eta_k) & =
  \setdef{\sigma \in \Sec{E}}{\sigma = \sigma_{B''}, B'' =
    \hspace*{-.25cm} \sum_{
    \begin{subarray}{c}
      \tilde{B} \in S(\psinv(B)) \\
      \cap\Br_{-1}(\famX) \cap \Br_{0}(\famX)
    \end{subarray}
  } \hspace*{-.75cm} (-1)^{\delta_{k+1}(\tilde{B})} \tilde{B} \, , \;
  B \in \Br^{2k}(\famX) \; \text{primitive}} \, .
\end{align*}
Let $\mathcal{C}_{\ver} (\eta;\eta_1,\dots,\eta_k) = \cup_{k \in
  \natural} \mathcal{C}^{(k)}_{\ver} (\eta;\eta_1,\dots,\eta_k) $,
$\mathcal{C}_{\hor}(\eta;\eta_1,\dots,\eta_k) = \cup_{k \in
  \natural} \mathcal{C}^{(k)}_{\hor} (\eta;\eta_1,\dots,\eta_k)$,
and denote by $C_{\ver}(\eta;\eta_1,\dots,\eta_k)$ and $C_{\hor}
(\eta;\eta_1,\dots,\eta_k)$, respectively, the subbundles of the Lie
algebroid $E$ generated by the latter families.

Taking into account the previous discussion, we are now ready to
compute~$\invclos{\{\Gamma -
  \eta\spV,\eta\spV_1,\dots,\eta\spV_k\}}$ for a mechanical control
system defined on a Lie algebroid.

\begin{proposition}
  \label{prop:decomposition-accessibility-distribution} 
  Let $m\in M$. Then,
  \begin{align*}
    & \invclos{\{\Gamma - \eta\spV,\eta\spV_1,\dots,\eta\spV_k\}} \cap
    \Ver_{0_m}(\prol{E}) = C_{\ver}(\eta;\eta_1,\dots,\eta_k)(m) \spV
    \, , \\
    & \invclos{\{\Gamma - \eta\spV,\eta\spV_1,\dots,\eta\spV_k \}}
    \cap \Hor_{m}(\prol{E}) = C_{\hor}(\eta;\eta_1,\dots,\eta_k)(m) \,
    .
  \end{align*}
\end{proposition}

\begin{proof}
  From the definition of the families $\mathcal{C}_{\ver}
  (\eta;\eta_1,\dots,\eta_k) $ and $\mathcal{C}_{\hor}
  (\eta;\eta_1,\dots,\eta_k) $, one sees that for each $B \in
  \Br(\famX)$ primitive, we have computed the $\real$-linear
  combinations from $\Br(\famX)$ that appear along with $B$ in the
  decomposition of $\psinv (B)$ according to~\eqref{eq:decomposition}.
  Since it is only these primitive brackets which appear when making
  Lie brackets of $\{\Gamma - \eta\spV_0,\eta\spV_1,\dots,\eta\spV_k
  \}$, this will generate $\invclos{\{\Gamma -
    \eta\spV_0,\eta\spV_1,\dots,\eta\spV_k \}}$ along the zero-section
  of $\prol{E}$.
\end{proof}

\begin{remark}
  {\rm In the absence of potential terms, i.e., $\eta =0$, one has
    that
    \begin{align*}
      C_{\ver}(0;\eta_1,\dots,\eta_k) = \symclos{\{
        \eta_1,\dots,\eta_k\}} \, , \quad
      C_{\hor}(0;\eta_1,\dots,\eta_k) =
      \invclos{\symclos{\{\eta_1,\dots,\eta_k\}}} \, .
    \end{align*}
    It is worth noticing that, in this case,
    $C_{\ver}(0;\eta_1,\dots,\eta_k) \subseteq
    C_{\hor}(0;\eta_1,\dots,\eta_k)$. This is not true in general.

}
\end{remark}

\subsection{Accessibility and controllability tests}

In this section we merge the notions introduced in
Section~\ref{se:notions} with the results obtained in
Section~\ref{se:structure} to give tests for accessibility and
controllability.

\begin{proposition}\label{prop:accessibility}
  Let $m \in M$ and assume the Lie algebroid $E$ is locally transitive
  at $m$.  Then the mechanical control sys\-tem~\eqref{eq:CS} is
  \begin{itemize}
  \item locally base accessible from $m$ if $\,
    C_{\hor}(\eta;\eta_1,\dots,\eta_k)(m) + \ker \rho = E_{m}$,
  \item locally accessible from $m$ at zero if $\,
    C_{\hor}(\eta;\eta_1,\dots,\eta_k)(m) + \ker \rho = E_{m}$ and $\,
    C_{\ver}(\eta;\eta_1,\dots,\eta_k)(m)= E_{m}$.
  \end{itemize}
\end{proposition}

\begin{proof}
  Consider the accessibility subbundle $\invclos{\{\Gamma -
    \eta\spV,\eta\spV_1,\dots,\eta\spV_k\}}$ in the Lie algebroid
  $\prol{E}$. Since $E$ is locally transitive at $m$ and $
  C_{\hor}(\eta;\eta_1,\dots,\eta_k)(m) + \ker \rho_m = E_{m}$ by
  hypothesis, and $C_{\hor}(\eta;\eta_1,\dots,\eta_k)(m) \subseteq
  \invclos{\{\Gamma - \eta \spV,\eta\spV_1,\dots,\eta\spV_k\}} (0_m)$
  (cf.
  Proposition~\ref{prop:decomposition-accessibility-distribution}),
  then $\rho^1 ( \invclos{\{\Gamma - \eta
    \spV,\eta\spV_1,\dots,\eta\spV_k\}} ) (0_m) = T_m(0_M(M))$. As a
  consequence, there exists an open connected submanifold $N$ of $M$
  containing $m$ such that $0_M(N)$ is an integral manifold of $\rho^1
  ( \invclos{\{\Gamma - \eta \spV,\eta\spV_1,\dots,\eta\spV_k\}} )$.
  Let $\lambda$ be the maximal integral manifold of $E$ which contains
  $0_M (N)$.  Since $\rho^1 (\invclos{\{\Gamma - \eta
    \spV,\eta\spV_1,\dots,\eta\spV_k\}} )$ is the accessibility
  distribution, $\lambda$ is invariant under~\eqref{eq:CS} and the
  system is locally accessible when restricted to $\lambda$. Thus the
  set $\R_E^U(0_{m},\le T)$ is open for all $U\subseteq \lambda$
  neighborhood of $0_{m}$ and sufficiently small $T$.  Let $V$ be a
  neighborhood of $m$ in $M$ and define $U = \tau^{-1}(V) \cap
  \lambda$. The result now follows from the fact that $\tau$ is an
  open mapping and hence the set $\tau (\R^U_E(0_{m},\le T)) \subset
  \R_M^V(m,\le T)$ is open in $M$ for $T$ sufficiently small.
  
  As for the second statement, note that the fact that $E$ is locally
  transitive at $m$ implies that $\prol{E}$ is locally transitive at
  $0_m$.  In addition, since $C_{\hor}(\eta;\eta_1,\dots,\eta_k)(m) +
  \ker \rho = E_{m}$ and $C_{\ver}(\eta;\eta_1,\dots,\eta_k)(m) =
  E_{m}$ by hypothesis, and $C_{\hor}(\eta;\eta_1,\dots,\eta_k)(m)$,
  $C_{\ver}(\eta;\eta_1,\dots,\eta_k)(m) \spV \subseteq
  \invclos{\{\Gamma - \eta \spV,\eta\spV_1,\dots,\eta\spV_k\}}
  (0_{m})$ by
  Proposition~\ref{prop:decomposition-accessibility-distribution}, we
  conclude $\invclos{\{\Gamma - \eta
    \spV,\eta\spV_1,\dots,\eta\spV_k\}} (0_{m}) + \ker \rho^1=
  \prol[0_m]{E}$.  The result now follows from
  Theorem~\ref{th:accessibility-test-in-E}.
\end{proof}

In the absence of potential energy terms, the accessibility tests
presented above simply read as follows.

\begin{proposition}\label{prop:accessibility-without-potential}
  Let $m \in M$ and assume the Lie algebroid $E$ is locally transitive
  at $m$.  A mechanical control sys\-tem~\eqref{eq:CS} with no
  potential terms is
  \begin{itemize}
  \item locally base accessible from $m$ if $\, \invclos{\symclos{\{
        \eta_1,\dots,\eta_k\}}}(m) + \ker \rho = E_{m}$.
  \item locally accessible from $m$ at zero if $\, \symclos{\{
      \eta_1,\dots,\eta_k\}} (m)= E_{m}$.
  \end{itemize}
\end{proposition}

\begin{remark}{\rm
    Alternatively, if the sufficient condition for locally base
    accessibility in
    Proposition~\ref{prop:accessibility-without-potential} is not met,
    i.e., $\invclos{\symclos{\{ \eta_1,\dots,\eta_k\}}}(m) + \ker \rho
    \not = E_{m}$, the corresponding proof also yields the following
    result.  Let $N$ denote the maximal integral manifold of
    $\invclos{\symclos{\{ \eta_1,\dots,\eta_k\}}}(m)$ passing through
    $m$.  Then, for each neighborhood $V$ of $m$ in $M$ and each $T$
    sufficiently small, $\R_M^V(m,\le T) \subset N$, and $\R_M^V(m,\le
    T)$ contains a non-empty open subset of~$N$.}
\end{remark}

The notions of good and bad symmetric products can be stated in a
similar way as for Lie brackets. We say that a symmetric product $P$
in the sections $\{ \eta,\eta_1,\dots,\eta_k \}$ is \emph{bad} if the
number of occurrences of each $\eta_i$ in $P$ is even. Otherwise, $P$
is \emph{good}. Accordingly, $\symprod{\eta_i}{\eta_i}$ is bad and
$\symprod{\symprod{\eta}{\eta_j}}{\symprod{\eta_i}{\eta_i}}$ is good.
The following theorem gives sufficient conditions for local
controllability.

\begin{proposition}\label{prop:controllability}
  Let $m \in M$.  The mechanical control sys\-tem~\eqref{eq:CS} is
  \begin{itemize}
  \item locally base controllable from $m$ if it is locally base
    accessible from $m$ and every bad symmetric product in $\{\eta,
    \eta_1,\dots,\eta_k \}$ evaluated at $m$ can be put as an
    $\real$-linear combination of good symmetric products of lower
    degree and elements of $\ker \rho$,
  \item locally controllable from $m$ at zero if it is locally
    accessible from $m$ at zero and every bad symmetric product in
    $\{\eta, \eta_1,\dots,\eta_k \}$ evaluated at $m$ can be put as an
    $\real$-linear combination of good symmetric products of lower
    degree.
  \end{itemize}
\end{proposition}

\begin{proof}
  The proof follows from the following considerations. First, note
  that every bad Lie bracket in $\{ \Gamma - \eta \spV
  \eta_1\spV,\dots,\eta_k\spV \}$ gives rise to the vertical lift of a
  bad symmetric product in $\{\eta,\eta_1,\dots,\eta_k \}$. Second,
  observe that every good symmetric product in $\{
  \eta,\eta_1,\dots,\eta_k \}$ can be alternatively written as a good
  Lie bracket in $\{ \Gamma - \eta \spV, \eta_1\spV,\dots,\eta_k\spV
  \}$ evaluated at the zero section $0_M$ (modulo a minus sign). The
  result is now an application of
  Proposition~\ref{th:controllability-test-in-E} to the setting of
  mechanical control systems on Lie algebroids.
\end{proof}

The corresponding tests for base accessibility and controllability
with regards to a manifold can be proved in a similar way.

\begin{proposition}\label{prop:accessibility-with-regards-to-N}
  Let $\psi:M \rightarrow N$ be an open map.  Let $m \in M$ and assume
  $\psi_* (\rho (E_m)) = T_{\psi(m)}N$.  Then the mechanical control
  sys\-tem~\eqref{eq:CS} is
  \begin{itemize}
  \item locally base accessible from $m$ with regards to $N$ if $\,
    C_{\hor}(\eta;\eta_1,\dots,\eta_k)(m) + \rho^{-1} (\ker \psi_*) =
    E_{m}$,
  \item locally base controllable from $m$ with regards to $N$ if the
    system is locally base accessible from $m$ with regards to $N$ and
    every bad symmetric product in $\{ \eta,\eta_1,\dots,\eta_k \}$
    evaluated at $m$ can be put as an $\real$-linear combination of
    good symmetric products of lower degree and elements of $\rho^{-1}
    (\ker \psi_*)$.
  \end{itemize}
\end{proposition}

\subsection{Mechanical control systems related by a morphism of Lie
  algebroids}\label{se:morphisms}

In this section we study the relation between the controllability
properties of two mechanical control systems related by a morphism of
Lie algebroids.  Consider a mechanical control system on $\tau:E
\rightarrow M$ with $\Gamma- \eta\spV$ as drift section (where
$\Gamma$ is a spray) and inputs $\{\eta\spV_1, \dots,\eta\spV_k\}$,
and a mechanical control system on $\tau:\ov{E} \rightarrow \ov{M}$
with $\ov{\Gamma}- \ov{\eta}\spV$ as drift section (where
$\ov{\Gamma}$ is a spray) and inputs
$\{\ov{\eta}\spV_1,\dots,\ov{\eta}\spV_{\ov{k}}\}$.

Let $\Psi:E \rightarrow \ov{E}$ be a morphism of Lie algebroids.
Assume the two mechanical control systems are weakly
$\prolmap{\Psi}$-related. Because of homogeneity, one can deduce that
$\prolmap{\Psi}\circ\Gamma = \ov{\Gamma} \circ \Psi$, so that $\Psi$
maps the corresponding associated connection $\nabla$ onto the
associated connection $\ov{\nabla}$. Moreover, using the definition of
morphism of Lie algebroids, one can conclude that
\begin{align*}
  \Psi \circ \eta = \ov{\eta} \circ \Psi + \sum_{j=1}^{\ov{k}} b^j
  (\ov{\eta}_j \circ \Psi) \, , \quad \Psi \circ \eta_i =
  \sum_{j=1}^{\ov{k}} C_i^j (\ov{\eta}_j \circ \Psi) \, , \quad i = 1
  , \dots, k \, ,
\end{align*}
for some functions $C_i^j$ on $M$, i.e., the relation by $\prol{\Psi}$
between the vertical lifts of the potential terms and the input
sections of the mechanical control systems translates into a relation
by $\Psi$ of the potential terms and the input sections themselves.
Following the steps of the proof of
Proposition~\ref{th:open-map-general}, one can also infer the
relationship between the base accessibility and controllability
properties of two $\Psi$-weakly related mechanical systems. As above,
the control functions of the second system are related to the control
functions of the first one by means of $\ov{u}_j(t)= b^j (m(t)) +
\sum_{i=1}^k C_i^j(m(t)) u_i(t)$.

\begin{proposition}\label{th:open-map}
  Let $\map{\Psi}{E}{\ov{E}}$ be a morphism of Lie algebroids such
  that the base map $\psi$ is an open map. Consider two mechanical
  control systems which are weakly $\Psi$-related. If the system on
  $E$ is locally base accessible (respectively locally base
  controllable) from $m$ then the system on $\ov{E}$ is locally base
  accessible (respectively locally base controllable) from $\psi(m)$.
\end{proposition}

In the particular case when $\Psi$ is an isomorphism between the
fibers of the Lie algebroids, if the two systems are $\Psi$-related,
then the sufficient conditions for local (base) accessibility
(respectively controllability) are either simultaneously satisfied on
$E$ and $\ov{E}$ or simultaneously not satisfied, as stated in the
following result.

\begin{proposition}\label{prop:equivalence}
  Let $\Psi:E \rightarrow \ov{E}$ be a morphism of Lie algebroids
  which is an isomorphism on each fiber. Consider two mechanical
  control systems on $E$ and $\ov{E}$, with $k \ge \ov{k}$, that are
  $\Psi$-related.  Let $m \in M$. Then
  \begin{enumerate}
  \item $C_{\ver}(\eta;\eta_1,\dots,\eta_k)(m) = E_m$ if and only if
    $C_{\ver}(\ov{\eta};\ov{\eta}_1,\dots,\ov{\eta}_{\ov{k}})
    (\psi(m)) = \ov{E}_{\psi(m)}$,
  \item $C_{\hor}(\eta;\eta_1,\dots,\eta_k)(m) + \ker \rho= E_m$ if
    and only if
    $C_{\hor}(\ov{\eta};\ov{\eta}_1,\dots,\ov{\eta}_{\ov{k}})
    (\psi(m)) + \ker \ov{\rho} = \ov{E}_{\psi(m)}$,
  \item Every bad symmetric product in $\{\eta, \eta_1, \dots, \eta_k
    \}$ evaluated at $m$ can be put as an $\real$-linear combination
    of good symmetric products of lower degree if and only if every
    bad symmetric product in $\{\ov{\eta}, \ov{\eta}_1, \dots,
    \ov{\eta}_{\ov{k}} \}$ evaluated at $\psi(m)$ can be put as an
    $\real$-linear combination of good symmetric products of lower
    degree.
  \end{enumerate}
\end{proposition}

\begin{proof}
  Note that the assumptions of the proposition imply $k= \ov{k}$.
  Therefore, the sections $\eta_i$ and $\ov{\eta}_i$ can be chosen to
  be $\Psi$-related
  \begin{equation}\label{eq:eta-reduced}
    \Psi \circ \eta_i = \ov{\eta}_i \circ \psi \, , \; 1 \le i \le k \, .
  \end{equation}
  Taking into account~\eqref{eq:eta-reduced}, one can verify that
  \begin{subequations}\label{eq:related-by-psi}
    \begin{align}
      & \Psi \circ C_{\ver}(\eta;\eta_1,\dots,\eta_k) =
      C_{\ver}(\ov{\eta}; \ov{\eta}_1, \dots, \ov{\eta}_{\ov{k}})
      \circ \psi \, ,
      \label{eq:related-by-psi-a} \\
      & \Psi \circ C_{\hor}(\eta;\eta_1,\dots,\eta_k) =
      C_{\hor}(\ov{\eta}; \ov{\eta}_1, \dots, \ov{\eta}_{\ov{k}})
      \circ \psi \, . \label{eq:related-by-psi-b}
    \end{align}
  \end{subequations}
  The proof now follows from~\eqref{eq:related-by-psi} and the
  hypothesis that $\Psi$ is an isomorphism on every fiber.
\end{proof}

\section{Applications to simple mechanical control systems and
  semidirect products}
\label{se:applications}
\setcounter{equation}{0}

In this section we show how the formalism of mechanical control
systems on Lie algebroids unifies the treatment of several situations
which have been previously considered in the literature. We recover
known accessibility and controllability results for simple mechanical
control systems and develop some new ones.

\subsubsection*{Simple mechanical control systems}

Let $Q$ be a $n$-dimensional manifold. A \emph{simple mechanical
  control system} is defined by a tuple
$(Q,\mathcal{G},V,\mathcal{F})$, where $Q$ is the manifold of
configurations of the system, $\mathcal{G}$ is a Riemannian metric on
$Q$ (the kinetic energy metric of the system), $V \in C^{\infty}(Q)$
is the potential function and $\mathcal{F}=\{F^1,\dots,F^k\}$ is a set
of $k$ linearly independent $1$-forms on $Q$, which physically
correspond to forces or torques.  The dynamics of simple mechanical
control systems is classically described by the forced
Euler-Lagrange's equations
\begin{equation}\label{eq:Lagrange}
  \frac{\partial }{\partial t} \left( \frac{\partial L}{\partial
  \dot{q}} \right) - \frac{\partial L}{\partial q} = \sum_{i=1}^k
  u_{i} (t) F^i \, , 
\end{equation}
where $L:TQ \rightarrow \real$, $L(q,\dot{q}) =\frac{1}{2} \mathcal{G}
(\dot{q},\dot{q}) - V(q)$ is the Lagrangian function of the system.

There are several ways of intrinsically writing these equations. Here,
we present a formulation following the Lie algebroid formalism
explained above.  Consider $E=TQ$, $M=Q$ and the mappings
$\tau=\tau_Q:TQ \rightarrow \real$, $\rho = \Id{TQ} :TQ \rightarrow
TQ$. Then, it is easy to see that $TQ$ is a Lie algebroid with anchor
map $\Id{TQ}$. In this setting, the forces in $\mathcal{F}$ correspond
to sections of the dual bundle $E^*=T^*Q$. By means of the musical
isomorphisms associated with the kinetic energy $\mathcal{G}$, we can
consider them as sections $\famY = \{ Y_1,\dots,Y_k \}$ of $E=TQ$
(i.e. vector fields). Then, equations~\eqref{eq:Lagrange} are
equivalently given by
\begin{equation}\label{eq:LCS-mechanics}
  \dot{a} (t) = \rho^1\Bigl(\Gamma (a(t)) - (\grad_{\Gc} V) \spV (a(t))+
  \sum_{i=1}^k u_i(t) Y_i \spV (a(t))\Bigr) \, , 
\end{equation}
where $\Gamma$ is the second order equation associated with
$\nabla^{\Gc}$. Here, the mapping $\rho^1$ is just the identity in
$TTM$.

The notion of base accessibility and controllability (resp.
accessibility and controllability at zero) in $M=Q$ precisely
corresponds to the concept of configuration accessibility and
controllability (resp. accessibility and controllability at zero
velocity) in $Q$ as introduced in~\cite{ADL-RMM:95c}.  Also, since
$\rho$ is an isomorphism, we conclude that the application of
Propositions~\ref{prop:accessibility} and~\ref{prop:controllability}
to this case just renders the known tests for accessibility and
controllability~\cite{ADL-RMM:95c}.

\subsubsection*{Simple mechanical control systems with symmetry I}

Assume that a simple mechanical control system
$(Q,\mathcal{G},V,\mathcal{F})$ is \emph{invariant under the action of
  a Lie group $G$ on~$Q$},
\[
\begin{array}{rrcl}
  \Phi:& G\times Q&\longrightarrow&Q \\
  & (g, q)&\longmapsto&\Phi(g,q)=\Phi_g(q)=gq \, .
\end{array}
\]
Invariance for the control system means that $\Phi^*_g \mathcal{G} =
\mathcal{G}$, $\Phi^*_g V = V$ and $\Phi^*_g F^i = F^i$, for $1 \le i
\le k$ and all $g \in G$. The orbit through a point $q$ is
$\hbox{Orb}_G(q)=\{gq\; |\; g\in G\}$. We denote by $\mathfrak{g}$ the
Lie algebra of $G$. For any element $\xi\in\mathfrak{g}$, let $\xi_Q$
denote the corresponding infinitesimal generator of the group action
on $Q$. Then,
\[
T_q(\hbox{Orb}_G(q))=\{\xi_Q (q)\; |\; \xi\in\mathfrak{g}\} \, .
\]
If the action $\Phi$ is free and proper, we can endow the quotient
space $Q/G$ with a manifold structure such that the canonical
projection $\pi: Q \longrightarrow Q/G$ is a surjective submersion.
Then, we have that $Q(Q/G, G, \pi)$ is a principal bundle with bundle
space $Q$, base space $Q/G$, structure group $G$ and projection $\pi$.
Note that the kernel of $T\pi$ consists of the vertical tangent
vectors, i.e., the vectors tangent to the orbits of $G$ in $Q$. We
denote the bundle of vertical vectors by $\mathcal{V}$, with
$\mathcal{V}_q=T_q (\hbox{Orb}_G(q))$, $q \in Q$.

The action $\Phi$ induces the \emph{lifted action of $G$ on $TQ$},
$\hat{\Phi}:G \times TQ \rightarrow TQ$, defined by
$\hat{\Phi}_g=T\Phi_g$.  Assuming $\Phi$ is free and proper, we have
that $\hat{\Phi}$ is also free and proper, and hence $p:TQ \rightarrow
TQ/G$, $p(v_q) = [v_q]$, is a surjective submersion. Consider then the
Lie algebroid defined by $E=TQ/G$, $M=Q/G$ and
\[
\tau ([v_q]) = [q] \, , \quad \rho([v_q]) = T\pi (v_q) \, .
\]
It is not difficult to verify that both $\tau$ and $\rho$ are
well-defined. The Lagrangian function $L$ being $G$-invariant, it
induces a reduced Lagrangian function $\ell :TQ/G \rightarrow \real$.
The invariance of the forces in $\mathcal{F}$, or equivalently the
fact that the sections in $\famY$ are invariant, implies that there
exist well-defined sections $\famB = \{ B_i:Q/G \rightarrow TQ/G \}$,
$1 \le i \le k$, such that $p \circ Y_i = B_i \circ \pi$. Finally, the
invariance of the potential function and the Riemannian metric implies
that there exists $\ov{\grad_{\Gc} V}$ such that $p \circ \grad_{\Gc}
V = \ov{\grad_{\Gc} V} \circ \pi$.  In this setting, the notion of
base accessibility and controllability precisely corresponds to
configuration accessibility and controllability in $Q/G$. However,
note that the notions of accessibility and controllability at zero in
$E$ are stronger than the notions of accessibility and controllability
at zero velocity in $Q/G$, since the former ones imply that the
reachable sets contain open sets in $E=TQ/G$, whereas the latter only
involve open sets in $T(Q/G)$.

If the mechanical control system is defined on a trivial principal
fiber bundle $Q=G \times Q/G$, one may consider the canonical
projection $\tau: G \times Q/G \rightarrow G$. This mapping is open,
and the notions of base accessibility and base controllability with
regards to $G$ (cf.  Section~\ref{se:notions}) precisely correspond to
the concepts of \emph{fiber configuration accessibility} and
\emph{fiber configuration controllability} as introduced
in~\cite{JC-SM-JPO-HZ:02}.
Proposition~\ref{prop:accessibility-with-regards-to-N} renders
appropriate tests to check these properties.

A different (and natural) question, however, concerns the search for
tests of accessibility or controllability in $Q$ which make use of the
symmetry properties of the mechanical control system. This is what we
analyze next.

\subsubsection*{Simple mechanical control systems with symmetry II}

Now, we will apply the results of Section~\ref{se:morphisms} to simple
mechanical control systems with symmetry. Denote by $E_1 \rightarrow
M_1$ the Lie algebroid $TQ \rightarrow Q$, and by $E_2 \rightarrow
M_2$ the Lie algebroid $TQ/G \rightarrow Q/G$. Note that both
algebroids have fibers of the same dimension $n = \dim Q$. Let
$\Psi:TQ \rightarrow TQ/G$ be the projection mapping associated with
the lifted action $\hat{\Phi}$, $\Psi = p$. It can be easily verified
that $\Psi$ defines a morphism of Lie algebroids.  The base mapping
$\psi:Q \rightarrow Q/G$ corresponds precisely to the projection
$\pi$. Finally, observe that $\Psi$ is surjective.

As a consequence of Proposition~\ref{prop:equivalence}(ii), we have
that the criterion to test base accessibility in $M_1=Q$ (cf.
Proposition~\ref{prop:accessibility}),
\[
C_{\hor}(\grad_{\Gc} V; \famY) = E_1 (=TQ) \, ,
\]
is equivalent to
\[
C_{\hor}(\ov{\grad_{\Gc} V}; \famB) = E_2 (=TQ/G) \, .
\]
Hence, in this way we simplify the computational cost to test the
controllability properties of the system, since one deals with the
reduced representation (and therefore, one works in a space of smaller
dimension). These are precisely the results obtained
in~\cite{JC-SM-JPO-HZ:02,SM:02} (although here the analysis is more general
since we consider nontrivial potential terms). The same simplification
occurs for accessibility at zero velocity (cf.
Proposition~\ref{prop:equivalence}(i)) ,
\[
C_{\ver}(\grad_{\Gc} V; \famY) = TQ \quad \text{if and only if} \quad
C_{\ver}(\ov{\grad_{\Gc} V}; \famB) = TQ/G \, .
\]
As for controllability in $Q$, Proposition~\ref{prop:equivalence}(iii)
ensures that it is enough to check that the bad symmetric products in
$\{\ov{\grad_{\Gc} V},B_1, \dots, B_k\}$ are $\real$-linear
combinations of good ones in $TQ/G$.  Finally,
Proposition~\ref{th:open-map} ensures that if the reduced system is
not base accessible (resp. controllable), then the original system is
not base accessible (resp. controllable).

\subsubsection*{Semidirect products}

Let $\G$ be a real Lie algebra acting transitively on a manifold $M$,
that is, let $\G\to\vectorfields{M}$ be a surjective Lie algebra
homomorphism mapping each element $\xi$ of $\G$ to a vector field
$\xi_M$ on $M$. Define the following Lie algebroid structure. The
bundle space is $E=M\times\G$ and the mapping $\tau: E \rightarrow M$
is just the projection onto the first factor.  The anchor map is given
by $\rho(m,\xi) = \xi_M(m)$.  The Lie bracket is defined via the
anchor map $\rho$ and the bracket of constant sections. The latter is
defined as the constant section corresponding to the bracket on $\G$,
that is, if $\sigma_1(m)=(m,\xi_1)$ and $\sigma_2(m)=(m,\xi_2)$ are
two constant sections, then
$\br{\sigma_1}{\sigma_2}(m)=(m,[\xi_1,\xi_2]_\G)$.  Note that the case
of a linear action on a vector space $\mathfrak{V}$ can be treated as
a semidirect product by taking an orbit of the action as the base
manifold $M$.

If we identify $TE\equiv TM\times T\G\equiv TM\times\G\times\G$ using
the left multiplication, an element of $\prol{E}$ is of the form
$(a,b,v)=\bigl((m,\xi),(m,\eta),(v_m,\xi,\zeta) \bigr)$. The condition
$T\tau(v) = \rho(b)$ simply implies that $v_m= \eta_M(m)$.  Therefore,
we can identify $\prol{E}$ with $M\times \G \times \G \times \G$, and
the corresponding maps are
\[
\tau_1(m,\xi,\eta,\zeta) = (m,\xi), \quad \prb(m,\xi,\eta,\zeta) =
(m,\eta), \quad \rho^1(m,\xi,\eta,\zeta) = (\eta_M(m),\xi,\zeta) \, .
\]

Let $(\Gc,V,\{\theta_1,\dots,\theta_k\})$ be a mechanical control
system on the Lie algebroid $E$. Assume the bundle metric $\Gc$ comes
from an inner product on $\G$ (and therefore does not depend on the
base point), $\Gc ((m,\xi_1),(m,\xi_2)) = \Gc (\xi_1,\xi_2)$. For each
$\xi \in \G$, define $\ada{\xi}: \G \rightarrow \G$ by $\Gc
(\ada{\xi}\eta_1,\eta_2) = \Gc (\eta_1,[\xi,\eta_2]_\G)$. The spray
associated with $\nabla^\Gc$ then reads
\[
\Gamma_{\nabla^\Gc} (m,\xi)=(m,\xi,\xi,\ada{\xi}\xi),
\]
and the controlled equations of motion are explicitly given by
\[
\dot{a} - \ada{a}{a} = -\grad_\Gc V(m)+\sum_{i=1}^k u_i\eta_i(m) \, .
\]
If one takes constant sections $\sigma_i(m)=(m,\xi_i)$, $i=1,2$, the
expression of the Levi-Civita connection is
\[
\nabla^\Gc_{\sigma_1}\sigma_2 (m) = \Bigl(m, \frac{1}{2}
[\xi_1,\xi_2]_{\mathfrak{g}} - \frac{1}{2}
(\ada{\xi_1}\xi_2+\ada{\xi_2}\xi_1) \Bigr) \,,
\]
and the symmetric product is
\[
\symprod{\sigma_1}{\sigma_2}(m) = \left( m,
  -(\ada{\xi_1}\xi_2+\ada{\xi_2}\xi_1) \right) .
\]

The above-developed tests can be easily applied to this kind of
problems in order to determine whether or not the system is base
accessible (resp. controllable).

Systems of the above type appears frequently as mechanical systems
defined on homogeneous spaces for a given group action. If a group $G$
acting transitively on a manifold $M$, one can consider the Lie
algebroid $TG\times M\rightarrow G\times M$, with anchor map
$\rho(v_g,m)=(v_g,0_m)$. The bracket of the Lie algebroid is just the
Lie bracket on the manifold $G$, where the coordinates of $M$ are
considered as \emph{parameters}.  Typically, one is thinking of
mechanical systems defined on $TG$ that depend on certain parameters
which are modeled by the coordinates of $M$. These systems are not
invariant under the right action on the group~$G$ on itself, but they
are invariant under the action on $G$ on itself and $M$ \emph{at the
  same time}. Therefore, one can consider the map $\map{\Psi}{TG\times
  M}{\G\times M}$, $\Psi(v_g,m)=(v_gg^{-1},gm)$, which is a morphism
of Lie algebroids with base map $\map{\psi}{G\times M}{M}$,
$\phi(g,m)=gm$.  One can verify that the map $\Psi$ is an isomorphism
in every fiber and that its associated base map $\psi$ is open.
Therefore, one can apply the results of Section~\ref{se:morphisms} to
determine whether a mechanical system on $TG$ depending on certain
parameters modeled by $M$ and which is not $G$-invariant, is locally
controllable by analyzing the corresponding system on the Lie
algebroid $\G \times M$.

An explicit example of the above type of system is the case of a rigid
body $G=SO(3)$ subject to control forces with fixed direction (in
space), which are the elements in $M$. Then the above procedure
precisely consists of studying the control problem in body
coordinates. The variables in $M$ evolve dynamically due to the
non-inertial rotating frame.

\section{Conclusions}\label{se:conclusions}

We have investigated the controllability properties of control systems
defined on Lie algebroids. We have established some general
controllability results for nonlinear affine control systems. We have
also introduced the concept of mechanical control system evolving on a
Lie algebroid. After defining appropriate accessibility and
controllability notions, we have investigated sufficient tests
guaranteeing them. We have paid special attention to the situation
where two control systems are related by means of a morphism of Lie
algebroids. Finally, we have illustrated the results with the class of
simple mechanical control systems and the class of systems evolving on
semidirect products. Future directions of research will include the
investigation of controllability tests along relative equilibria of
mechanical control systems on Lie algebroids and the treatment of
models that include gyroscopic forces and dissipation.

\section*{Acknowledgments}
The first author wishes to thank Sonia Mart{\'\i}nez for helpful
conversations and constant support. J. Cort\'es' work was partially
funded by NSF grant CMS-0100162 and EU Training and Mobility of
Researchers Program ERB FMRXCT-970137. E. Mart{\'\i}nez's work was
partially funded by CICYT grant BFM2000-1066-C03-01.

\section{Appendix}\label{se:appendix}

Here we gather some basic definitions concerning control systems
defined on manifolds~\cite{HN-AJvdS:90}. Let $\R_M^V (m,T)$ be the
reachable set from a point $m \in M$ at time $T >0$, following
trajectories which remain in the neighborhood $V$ of $m$ in $M$ for $t
\le T$.  Denote
\[
\R_M^V (m,\le T) = \bigcup_{t \le T} \R_M^V (m, t) \, .
\]
  
\begin{definition}
  The system~\eqref{eq:control-system} is \emph{locally accessible
    from $m \in M$} if $\R_M^V(m,\le T)$ contains a non-empty open set
  of $M$ for all neighborhoods $V$ of $m$ and all $T>0$. If this holds
  for any $m \in M$, then the system is called \emph{locally
    accessible}.
\end{definition}

\begin{definition}
  The system~\eqref{eq:control-system} is \emph{locally controllable
    from $m \in M$} if $\R_M^V(m,\le T)$ contains a non-empty open set
  of $M$ to which $m$ belongs for all neighborhoods $V$ of $m$ and all
  $T>0$.  If this holds for any $m \in M$, then the system is called
  \emph{locally controllable}.
\end{definition}
  
The \emph{accessibility algebra $\C$ of the control
  system~\eqref{eq:standard}} is defined as the smallest subalgebra of
$\vectorfields{M}$ containing $f,g_1,\dots,g_k$. It is not difficult to show
that every element of $\C$ is a linear combination of repeated Lie brackets
of the form
\[
[X_l,[X_{l-1},[\dots,[X_2,X_1]\dots]]] \, ,
\]
where $X_i \in \{f,g_1,\dots,g_k\}$, $1 \le i \le l$ and $l \in \natural$.\
The \emph{accessibility distribution} $C$ is defined as the distribution on
$M$ generated by the accessibility algebra $\C$,
\[
C(m) = \spn \left\{ X(m) \; | \; X \; \hbox{vector field in} \; \C
\right\} , \quad m \in M \, .
\]
\end{document}